\documentclass{article}
\usepackage{graphicx,bm,url,amssymb,amsmath,hyperref,color}
\usepackage{geometry}

\def\no{\noindent}
\def\pmatrix{\left(\begin{array}}
\def\endpmatrix{\end{array}\right)}

\def\CC{\mathbb{C}}

\def\NN{\mathbb{N}}

\def\RR{\mathbb{R}}

\def\D{{\cal D}}

\def\I{{\cal I}}

\def\P{{\cal P}}

\def\dd{\mathrm{d}}

\def\diag{\mathrm{diag}}

\def\rank{\mathrm{rank}}

\def\ii{\mathrm{i}}

\def\arg{\mathrm{arg}}

\newtheorem{theo}{Theorem}
\newtheorem{lem}{Lemma}

\newtheorem{rem}{Remark}
\newtheorem{defi}{Definition}
\def\proof{\noindent\underline{Proof}\quad}
\def\QED{\mbox{$\hfill{}~ \Box$}}

\def\bfb{{\bm{b}}}
\def\bfc{{\bm{c}}}

\def\bfgamma{{\bm{\gamma}}}

\def\eps{\varepsilon}

\def\aa{{\alpha}}
\def\lam{{\lambda}}

\def\cpr{$^\copyright\,$}

\begin{document}

\title{Analysis and implementation of collocation methods for fractional differential equations}

\author{
Luigi Brugnano\,\footnote{Dipartimento di Matematica e Informatica ``U.\,Dini'', 
             Universit\`a di Firenze,  Italy,~    \url{luigi.brugnano@unifi.it}} 
   \and  
Gianmarco Gurioli\,\footnote{Dipartimento di Matematica e Informatica ``U.\,Dini'',
             Universit\`a di Firenze,  Italy, ~ \url{gianmarco.gurioli@unifi.it}} 
     \and 
Felice Iavernaro\,\footnote{Dipartimento di Matematica, 
           Universit\`a di Bari Aldo Moro, Italy, ~ \url{felice.iavernaro@uniba.it}} 
     \and  
Mikk Vikerpuur\,\footnote{Institute of Mathematics and Statistics, 
           University of Tartu, Estonia, ~ \url{mikk.vikerpuur@ut.ee}} 
            }

\maketitle

\begin{abstract} 
Recently, the class of Runge-Kutta type methods named Fractional HBVMs (FHBVMs) has been introduced for the numerical solution of initial value problems of fractional differential equations, and a corresponding Matlab\cpr software has been released. Though an error analysis has already been given, a corresponding linear stability analysis is still lacking. We here provide such an analysis, together with some improvements concerning the mesh selection. This latter has been  implemented into a new version of the code, which is available on the web.

\medskip
\no{\bf Keywords:} Fractional Differential Equations, FDEs, Caputo derivative, Fractional HBVMs, FHBVMs.

\medskip
\no{\bf MSC:}  34A08, 65R20, 65-04.
\end{abstract}

\section{Introduction}

The mathematical modeling of Fractional Differential Equations (FDEs) has gained more and more importance in a number of applications in several scientific settings: we refer, e.g., to \cite{BegOrs05,Ben00,BueKay14,Bue20,BurOro14,BurOr17,Caf12,Cap13,Cus15,Garra11,Hen08,Hori12,Koh90,Magi10,NigArb07,Ors09,SaiZas97,She11,Ucha22} to mention a few of them (see also the review papers \cite{MKM2011,SZBCC2018}), and to the classical monographs \cite{Diethelm2010,Pod99} for an introduction to the subject.

Since finding analytical solutions of FDEs can be even more challenging than solving standard ordinary differential equations, there has been a growing need for suitable numerical methods, aiming to go beyond a first order approximation to the Caputo derivative (see \cite{Garr18} for a review on computational methods for FDEs). A first effort to consider higher-order approximations for the Caputo derivative was done in \cite{Lub85}, achieving convergence rates of the order of the underlying multistep method, also in the generic case of solutions that are not smooth at the initial time. From the computational viewpoint, solving fractional differential equations in an accurate, reliable, and efficient way can be arduous, because of the nonlocality of the operator. Moreover, since each time step requires evaluating contributions from all previous steps, making long-time simulations may be computationally demanding. Consequently,  the efficient treatment of the persistent memory term is an issue in itself, together with the solution of the nonlinear systems involved in implicit methods. To mitigate the former problem, approaches like fixed time windowing \cite{Pod99} and efficient discrete convolution techniques \cite{Sha06,Zen18} have been proposed, along with memory-saving strategies that split the fractional operator into a local part with fixed memory and a history part computed via Laguerre-Gauss quadrature, with a stepsize-independent error \cite{Zeng18}. In this context, it is worth noting that constant stepsizes may be inefficient, since the solutions of FDEs may possess a singularity at the initial time instant, leading to the need of correction terms \cite{Lub86}, the use of a smoothing transformation \cite{FPV2023,PTV2017}, or the adoption of graded meshes \cite{Li16,Sty17,Zen17,Yus12}  in place of uniform ones. Generalizations of the classical one-step Adams-Bashforth/Moulton scheme for first-order equations \cite{Die04,Die05}, scalable techniques and finite element/difference methods \cite{Bur12,Li16,Yus05}, trapezoidal methods \cite{Garr11}, product integration rules \cite{GarPop11}, and Krylov \cite{Mor11} methods have been also investigated in the literature, as well as collocation-type methods \cite{LS2018,PT2014}.

More recently, the class of {\em Fractional HBVMs (FHBVMs)} \cite{BBBI2024} has been introduced for the efficient numerical solution of initial value problems of fractional differential equations (FDE-IVPs) in the form (hereafter, we shall consider $ \aa\in(\ell-1,\ell)$, where $\ell\in\NN$, $\ell\ge1$, and $y\in A^\ell([0,T])$, i.e., $y^{(\ell)}$ is absolutely continuous) 
\begin{equation}\label{fde}
y^{(\aa)}(t) = f(y(t)), \qquad t\in[0,T], \qquad y^{(j)}(0) = y_0^j\in\RR^m,\quad j=0,\dots,\ell-1,
\end{equation}
with 
$$y^{(\aa)}(t) = \frac{1}{\Gamma(\ell-\aa)}\int_0^t(t-x)^{\ell-1-\aa}y^{(\ell)}(x)\dd x$$  
the Caputo fractional derivative of $y$.\footnote{For the sake of brevity, and withouth loss of generality, we have omitted $t$ as a formal argument of the vector field $f$.} In such a case, it is known that the solution of (\ref{fde}) is given by:
\begin{equation}\label{solfde}
y(t) ~=~ T_\ell(t) + I^\aa f(y(t)) ~\equiv~ \sum_{j=0}^{\ell-1}\frac{t^j}{j!}y_0^j ~+~ \frac{1}{\Gamma(\aa)}\int_0^t (t-x)^{\aa-1}f(y(x))\dd x,
\qquad t\in[0,T],
\end{equation}
where $I^\aa f(y(t))$ is the corresponding Riemann-Liouville integral. The methodology introduced in \cite{BBBI2024} is aimed at obtaining spectrally accurate solutions, following the results in the ODE case \cite{ABI2020}. The approach consists in the expansion of the vector field along the Jacobi polynomial basis $\{P_j\}_{j\ge0}$, orthonormal on the interval $[0,1]$ w.r.t. the weighting function 
\begin{equation}\label{orto}
\omega(c) = \aa(1-c)^{\aa-1}, \quad c\in[0,1], \quad \Longrightarrow \quad \int_0^1 \omega(c) P_i(c)P_j(c)\dd c = \delta_{ij},\quad i,j,=0,1,\dots.
\end{equation}
Considering a polynomial approximation of degree $s-1$ for the vector field, and approximating the corresponding Fourier coefficients by using the interpolatory Gauss-Jacobi quadrature of order $2k$ based at the zeros of $P_k$, one eventually obtains a FHBVM$(k,s)$ method \cite{BBBI2024}. All implementation details are given in \cite{BGI2024} (see also \cite{BGI2025}), and the Matlab\cpr code {\tt fhbvm} has been made available at the URL \cite{fhbvm}. In particular, the code {\tt fhbvm} automatically chooses the mesh to be used, i.e.:
\begin{itemize}
\item either a uniform mesh with timestep
\begin{equation}\label{h} h = \frac{T}N,\end{equation}
when the vector field is smooth at the origin;
\item or a graded mesh with the increasing timesteps
\begin{equation}\label{hi}
h_i = r^{i-1} h_1, \qquad i=1,\dots,\nu,
\end{equation}
for a suitable $r>1$, when the vector field is not smooth at the origin.
\end{itemize}
According to the analysis in \cite{BBBI2024}, the global error is bounded by $O(h^{s+\aa-1})$, in the case (\ref{h}), or $O(h_1^{2\aa}+h_\nu^{s+\aa})$, in the case (\ref{hi}). This choice, though generally very effective, can be relatively inefficient in case where $T$ in (\ref{fde}) is very large, and the solution of the given problem is both oscillatory  and nonsmooth at the origin. In fact, in such a case, an initial graded mesh, coupled with a subsequent uniform one, would be much more appropriate. Consequently, we shall here consider the case where the initial graded mesh is given by (\ref{hi}), with the constraint, for a suitable $n\in\{1,\dots,N\}$:
\begin{equation}\label{rnu}
\sum_{i=1}^\nu h_i ~\equiv~ h_1\sum_{i=1}^\nu r^{i-1} ~=~ h_1\frac{r^\nu-1}{r-1} ~=~ n h, 
\end{equation}
being $h$ the constant timestep defined in (\ref{h}), later used in the interval $[nh,T]$. In other words, we use a graded mesh only in the interval $[0,n h]$: after that, a uniform mesh with timestep $h$ is considered. This strategy apparently increases the complexity of the pre-processing phase of the code, as we shall see later, but overall it will confer a much greater efficiency. 

With this premise, in Section~\ref{fhbvm_sec} we recall the main facts about FHBVMs and, more specifically, the particular case where they reduce to collocation methods. In Section~\ref{linsys} we sketch a linear stability analysis of the methods, which is still missing. After that, in Section~\ref{appros} we provide full details for the practical implementation of the mixed stepsize strategy. Subsequently, in Section~\ref{discre} we sketch the discretization procedure implemented in the new Matlab\cpr code {\tt fhbvm2}. Some numerical tests, showing the efficiency of the new code, are then reported in Section~\ref{numtest}. At last, some conclusions are given in Section~\ref{fine}. 

\section{Fractional HBVMs (FHBVMs)}\label{fhbvm_sec}
In order to briefly introduce FHBVMs, let us consider the solution of problem (\ref{fde}) over the interval $[0,h]$, which we can rewrite as
\begin{equation}\label{fde1}
y^{(\aa)}(ch) = f(y(ch)), \qquad c\in[0,1],\qquad y^{(j)}(0) = y_0^j\in\RR^m,\quad j=0,\dots,\ell-1.
\end{equation}
Considering the expansion of the vector field along the orthonormal polynomial basis (\ref{orto}), one then obtains that (\ref{fde1}) can be rewritten as:
\begin{equation}\label{fde2}
y^{(\aa)}(ch) = \sum_{j\ge0} P_j(c)\gamma_j(y), \qquad c\in[0,1],\qquad y^{(j)}(0) = y_0^j\in\RR^m,\quad j=0,\dots,\ell-1,
\end{equation}
with the Fourier coefficients defined by
\begin{equation}\label{gamj}
\gamma_j(y) = \int_0^1 \omega(\tau)P_j(\tau)f(y(\tau h))\dd\tau, \qquad j=0,1,\dots.
\end{equation}
Integrating side by side, and imposing the initial conditions, one then obtains
$$
y(ch) = T_\ell(ch) + h^\aa\sum_{j\ge0} I^\aa P_j(c)\gamma_j(y), \qquad c\in[0,1],
$$
which is equivalent to (\ref{solfde}). In particular, by considering that 
$$P_0(c)\equiv 1 \qquad \mbox{and}\qquad I^\aa P_j(1)=\frac{\delta_{j0}}{\Gamma(\aa+1)},$$ 
one derives:
\begin{eqnarray}\nonumber
y(h) &=& T_\ell(h) + \frac{h^\aa}{\Gamma(\aa+1)}\gamma_0(y) ~=~T_\ell(h) + \frac{h^\aa}{\Gamma(\aa)\aa}\gamma_0(y)
\\[1mm] \label{solh}
&=& T_\ell(h) + \frac{1}{\Gamma(\aa)}\int_0^h (h-x)^{\aa-1}f(y(x))\dd x ~\equiv~ T_\ell(h) + I^\aa f(y(h)),
\end{eqnarray}
that is, (\ref{solfde}) at $t=h$. 
A polynomial approximation of degree $s-1$ to (\ref{fde2}) is obtained by truncating the infinite series in (\ref{fde2}) to a finite sum with $s$ terms:
\begin{equation}\label{sigalfa}
\sigma^{(\aa)}(ch) = \sum_{j=0}^{s-1} P_j(c)\gamma_j(\sigma), \qquad c\in[0,1],\qquad \sigma^{(j)}(0) = y_0^j\in\RR^m,\quad j=0,\dots,\ell-1,
\end{equation}
with the Fourier coefficients $\gamma_j(\sigma)$ defined according to (\ref{gamj}) by formally replacing $y$ with $\sigma$. As is clear, this is equivalent to requiring that the residual be orthogonal to all polynomials of degree $s-1$, w.r.t. the inner product defined by (\ref{orto}):\footnote{This feature is common to HBVMs, obtained when $\aa=1$ \cite{ABI2023,book,BI2022}.}
$$\int_0^1\omega(c)P_j(c)\left[ \sigma^{(\aa)}(ch)-f(\sigma(ch))\right]\dd c = 0, \qquad j=0,\dots,s-1.$$
The approximation $\sigma$, in turn, is given by integrating (\ref{sigalfa}) side by side and imposing the initial conditions:
\begin{equation}\label{sig}
\sigma(ch) = T_\ell(ch) + h^\aa\sum_{j=0}^{s-1} I^\aa P_j(c)\gamma_j(\sigma), \qquad c\in[0,1].
\end{equation}
In particular, at $t=h$ one derives, similarly as in (\ref{solh}):
\begin{eqnarray*}
\sigma(h) &=& T_\ell(h) +\frac{h^\aa}{\Gamma(\aa+1)}\gamma_0(\sigma) \\ 
&=& T_\ell(h) + \frac{1}{\Gamma(\aa)}\int_0^h (h-x)^{\aa-1}f(\sigma(x))\dd x ~\equiv~ T_\ell(h) + I^\aa f(\sigma(h)).
\end{eqnarray*} 
However, in order to derive a practical numerical method, the Fourier coefficients $\gamma_j(\sigma)$ need to be approximated by using a suitable quadrature rule. As anticipated in the introduction, by considering, for this purpose, for a suitable $k\ge s$,   the interpolatory Gauss-Jacobi quadrature of order $2k$, with abscissae $c_1,\dots,c_k$, placed at the zeros of $P_k(c)$,
$$P_k(c_i) = 0, \qquad i=1,\dots,k,$$
and corresponding weights $b_1,\dots,b_k$, one derives a FHBM$(k,s)$ method \cite{BBBI2024}. 
Such a method, in turn, admits a Runge-Kutta type formulation. In fact, by setting\,\footnote{For sake of brevity, we shall continue using the same symbol $\sigma$ both for denoting the approximation (\ref{sig}) and the one after the discretization of the integrals for computing the Fourier coefficients.} $$Y_i :=\sigma(c_ih), \qquad i=1,\dots,k,$$ one obtains that
\begin{equation}\label{hgamj}
\gamma_j(\sigma) \,=\,\int_0^1 \omega(c) P_j(c)f(\sigma(c h))\dd c
~\approx~ \sum_{\nu=1}^k b_\nu P_j(c_\nu) f(Y_\nu) \,=:\, \hat\gamma_j, \qquad j=0,\dots,s-1.
\end{equation}
Consequently, from (\ref{sig}), one derives the following stage equations:
\begin{eqnarray}\nonumber
Y_i &=& T_\ell(c_ih) + h^\aa\sum_{j=0}^{s-1}I^\aa P_j(c_i) \sum_{\nu=1}^k b_\nu P_j(c_\nu) f(Y_\nu)\\ \label{Yi}
&=& T_\ell(c_ih) + h^\aa\sum_{j=1}^k b_j \left[ \sum_{\nu=0}^{s-1} I^\aa P_\nu(c_i) P_\nu(c_j)\right] f(Y_j), \qquad i=1,\dots,k,
\end{eqnarray}
with the new approximation given by
\begin{equation}\label{y1}
y_1 :=\sigma(h) = T_\ell(h) + \frac{h^\aa}{\Gamma(\aa+1)}\sum_{i=1}^k b_i f(Y_i).
\end{equation}
The following result can be proved \cite{BBBI2024}.

\begin{theo}\label{RKform}
Let $\bfc=(c_1,\dots,c_k)^\top$, $\bfb=(b_1,\dots,b_k)^\top$ be the vectors with the abscissae and weights of the quadrature, respectively, and, moreover, let us set $$Y=\pmatrix{c}Y_1\\ \vdots\\ Y_k\endpmatrix,~ f(Y) = \pmatrix{c}f(Y_1)\\ \vdots\\ f(Y_k)\endpmatrix, ~T_\ell(\bfc h) = \pmatrix{c} T_\ell(c_1h)\\ \vdots\\ T_\ell(c_kh)\endpmatrix\in\RR^{km}.$$ Then, the equations (\ref{Yi})-(\ref{y1}) can be rewritten, respectively, as:\footnote{Hereafter, when not differently stated, $I$ will denote the identity matrix having the same size $m$ of the continuous problem (\ref{fde}). Differently, $I_r$ will denote the identity matrix of dimension $r$.}
\begin{equation}\label{RKeq}
Y = T_\ell(\bfc h) + h^\aa \I_s^\aa\P_s^\top\Omega\otimes I f(Y),\qquad y_1 = T_\ell(h) + \frac{h^\aa}{\Gamma(\aa+1)}\bfb^\top\otimes I f(Y),
\end{equation}
where $\Omega = \diag(\bfb)\in\RR^{k\times k}$, and
\begin{equation}\label{IP}
\I_s^\aa = \pmatrix{ccc} I^\aa P_0(c_1), &\dots& I^\aa P_{s-1}(c_1)\\
\vdots & &\vdots\\ I^\aa P_0(c_k), &\dots& I^\aa P_{s-1}(c_k)\endpmatrix,~
\P_s = \pmatrix{ccc} P_0(c_1), &\dots& P_{s-1}(c_1)\\
\vdots & &\vdots\\ P_0(c_k), &\dots& P_{s-1}(c_k)\endpmatrix\in\RR^{k\times s}.
\end{equation}
\end{theo}
\bigskip

\begin{rem} From the statement of Theorem~\ref{RKform}, one has that the Butcher tableau of a FHBVM$(k,s)$ method is given by 
\begin{equation}\label{Butab}
\begin{array}{c|c} \bfc & A\\ \hline \\[-3mm] &\bfb^\top\end{array}, \qquad A ~:=~ \I_s^\aa \P_s^\top\Omega.
\end{equation}
In particular, when $\aa=1$, it reduces to that of a HBVM$(k,s)$ method \cite{book}. Moreover, it can be proved that
\begin{equation}\label{rankA}
\forall k\ge s:~ \rank(A) ~ = ~ s.
\end{equation}
\end{rem}

It is worth mentioning that the discrete problem (\ref{RKeq}), which has (block) dimension $k$, can always be recasted into an equivalent problem of (block) dimension $s$, {\em independently} of $k$.\footnote{This feature is inherited from HBVMs \cite{book}, obtained in the case where $\aa=1$.} As a matter of fact, by setting
\begin{equation}\label{bgam}
\hat\bfgamma = \pmatrix{c}\hat\gamma_0\\ \vdots\\ \hat\gamma_{s-1}\endpmatrix \in\RR^{sm}
\end{equation}
the block vector with the approximate $s$ Fourier coefficients, from (\ref{hgamj}) it follows that
\begin{equation}\label{bgam1}
\hat\bfgamma = \P_s^\top\Omega\otimes I f(Y).
\end{equation}
As a result, the equations (\ref{RKeq}) can be recast, respectively, as
\begin{equation}\label{RKeq1}
Y = T_\ell(\bfc h) + h^\aa\I_s^\aa\otimes I \,\hat\bfgamma, \qquad y_1 = T_\ell(h) + \frac{h^\aa}{\Gamma(\aa+1)}\hat\gamma_0.
\end{equation}
Consequently, by plugging the second member of the first equation in (\ref{RKeq1}) to the second member in (\ref{bgam1}), one derives the equivalent discrete problem
\begin{equation}\label{bgam2}
\hat\bfgamma = \P_s^\top\Omega\otimes I\, f\left(T_\ell(\bfc h) + h^\aa\I_s^\aa\otimes I \,\hat\bfgamma\right),
\end{equation}
having (block) dimension $s$. Once it is solved, the new approximation is obtained by the second equation in (\ref{RKeq1}). We refer to \cite{BGI2024} for an efficient procedure for numerically solving (\ref{bgam2}).

\subsection{Collocation methods}\label{collocation}

To begin with, let us state this preliminary result.

\begin{lem}\label{POlem}
With reference to the matrices $\P_s$ and $\Omega$ defined in Theorem~\ref{RKform}, one has:
$$\P_s^\top\Omega \P_s = I_s\in\RR^{s\times s}.$$
Consequently, \quad $k=s \quad \Rightarrow\quad \P_s^\top\Omega = \P_s^{-1}$.
\end{lem}
\proof By taking into account that the quadrature rule has order $2k$, and therefore is exact for polynomial integrands of degree $2k-1$, one has, by setting $e_i,e_j\in\RR^s$ the $i$th and $j$th unit vectors, respectively, and taking into account (\ref{orto}) and (\ref{IP}):
\begin{eqnarray*}
e_i^\top\left(\P_s^\top\Omega \P_s\right)e_j &=& (\P_s e_i)^\top \Omega (\P_s e_j) ~=~ \sum_{\nu=1}^k b_\nu P_{i-1}(c_\nu)P_{j-1}(c_\nu) \\
&=&\int_0^1 \omega(c)P_{i-1}(c)P_{j-1}(c)\dd c ~=~\delta_{ij}, \qquad \forall i,j=1,\dots,s.
\end{eqnarray*}
Consequently, the first part of the statement follows. The statement is completed by observing that, when $k=s$, matrix $\P_s\in\RR^{s\times s}$.\,\QED\bigskip

Let us now denote
\begin{equation}\label{sigmalfa}
\sigma^{(\aa)}(\bfc h) = \pmatrix{c} \sigma^{(\aa)}(c_1h)\\ \vdots\\ \sigma^{(\aa)}(c_kh)\endpmatrix, \qquad
\sigma(\bfc h) = \pmatrix{c} \sigma(c_1h)\\ \vdots\\ \sigma(c_kh)\endpmatrix.
\end{equation}
As is clear, $\sigma(\bfc h)=Y$, the stage vector of the method. The following result then holds true.

\begin{lem}\label{sigalf}
With reference to (\ref{IP}) and (\ref{sigmalfa}), one derives
$$\sigma^{(\aa)}(\bfc h) = \P_s\P_s^\top\Omega\otimes I\,f(\sigma(\bfc h)).$$
\end{lem}
\proof The statement easily follows by evaluating the polynomial approximation,
$$\sigma^{(\aa)}(ch) = \sum_{j=0}^{s-1} P_j(c)\hat\gamma_j, \qquad c\in[0,1],$$
at the abscissae $c_1,\dots,c_k$, then considering (\ref{bgam})-(\ref{bgam1}).\,\QED\bigskip

We can now state the following result.

\begin{theo}\label{coll}
The FHBVM$(s,s)$ method collocates problem (\ref{fde}) at the abscissae $c_1h,\dots,c_sh$.
\end{theo}
\proof In fact, since $k=s$, from Lemmas~\ref{POlem} and \ref{sigalf} one has:
$$\sigma^{(\aa)}(\bfc h) = \P_s\P_s^\top\Omega\otimes I\,f(\sigma(\bfc h)) = \P_s\P_s^{-1}\otimes I\,f(\sigma(\bfc h))=f(\sigma(\bfc h)),$$
i.e. $$\sigma^{(\aa)}(c_ih) = f(\sigma(c_ih)), \qquad i=1,\dots,s.$$
As is clear, $\sigma$ also satisfies the initial conditions in (\ref{fde}) and, moreover, $y_1:=\sigma(h)$.\,\QED

\medskip
\begin{rem} It is worth noticing that, when $k=s$, the classical stage equation (\ref{RKeq}) and the equivalent formulation (\ref{bgam2}) of the discrete problem, do have the same (block) dimension.\end{rem}

\section{Linear stability analysis}\label{linsys}

As in the ODE case, also for FDEs it is customary to study the behavior of numerical methods when solving the linear test equation
\begin{equation}\label{test}
y^{(\aa)} = \lam y, \qquad y(0) = y_0\in\RR,\qquad \aa\in(0,1), \qquad \lam\in\CC, 
\end{equation}
whose solution is given by
\begin{equation}\label{sol}
y(t) = E_\aa(\lam t^\aa)y_0, \qquad t\ge0,
\end{equation}
where 
\begin{equation}\label{Ea}
E_\aa(t) = \sum_{j\ge0} \frac{t^j}{\Gamma(\aa j+1)}
\end{equation}
is the one-parameter Mittag-Leffler function.\footnote{As is clear, when $\aa=1$, from (\ref{Ea})  one retrieves the exponential function.}  This is a special case of the more general two-parameters Mittag-Leffler function,
\begin{equation}\label{Eab}
E_{\aa,\beta}(t) = \sum_{j\ge0} \frac{t^j}{\Gamma(\aa j+\beta)},\qquad \beta>0,
\end{equation}
since $E_\aa(t) \equiv E_{\aa,1}(t)$. The following properties are known (\cite{M1996,WZ2023}):
\begin{equation}\label{Ea0}
E_\aa(\lam t^\aa) \rightarrow 0, \mbox{~as~} t\rightarrow \infty \qquad \Longleftrightarrow\qquad
\lam\in\Lambda_\aa = \left\{ z\in\CC\backslash\{0\}\,:\,|\arg(z)|>\aa\frac{\pi}2\right\}.
\end{equation}
In particular,
\begin{equation}\label{talfa}
\lam\in\Lambda_\aa \qquad \Longrightarrow\qquad E_\aa(\lam t^\aa) = O(t^{-\aa}), \mbox{\quad as\quad }t\rightarrow\infty.
\end{equation}
Moreover, a simple calculation shows that (see (\ref{Eab}))
\begin{equation}\label{dEa}
\frac{\dd}{\dd t} E_\aa(\lam t^\aa) = \lam t^{\aa-1}E_{\aa\aa}(\lam t^\aa).
\end{equation}
From (\ref{talfa}) and (\ref{dEa}), one derives
\begin{equation}\label{talfa1}
\lam\in\Lambda_\aa \qquad \Longrightarrow\qquad t^{\aa-1}E_{\aa\aa}(\lam t^\aa) = O(t^{-\aa-1}), \mbox{\quad as\quad}t\rightarrow\infty,
\end{equation}
and, consequently,
\begin{equation}\label{roa}
\rho_\aa^\lam := \int_0^\infty t^{\aa-1}|E_{\aa\aa}(\lam t^\aa)|\dd t<\infty.
\end{equation}

\subsection{The continuous case}\label{continuous}

Stability results for the test equation (\ref{test}) allow to discuss the more general equation
\begin{equation}\label{testf}
y^{(\aa)} = \lam y + f(y), \qquad y(0) = y_0,
\end{equation}
where $f(y)$ is assumed to satisfy $f(0)=0$ (and, consequently, the origin is an equilibrium for (\ref{testf})). Moreover, we assume $f(y)$ to be Lipschitz with constant $L$. By considering that, by the nonlinear variation of constants formula, the solution of (\ref{testf}) is given by
\begin{equation}\label{solf}
y(t) = E_\aa(\lam t^\aa)y_0 + \int_0^t (t-s)^{\aa-1}E_{\aa\aa}(\lam(t-s)^\aa)f(y(s))\dd s, \qquad t\ge0, 
\end{equation}
the following stability result  by first approximation can be proved, by suitably adapting the arguments in \cite{CDT2018} (also, see that reference for possible alternative formulations).

\begin{theo} Assume that the origin is asymptotically stable for the test equation (\ref{test}) (i.e., $\lam\in\Lambda_\aa$) and, moreover, with reference to (\ref{roa}),\footnote{Actually, $L\rho_\aa^\lam<1$ would be enough.} 
\begin{equation}\label{minor}
L\rho_\aa^\lam < \frac{1}2. 
\end{equation}
Then, the origin is asymptotically stable for the complete equation (\ref{testf}).
\end{theo}
\proof 
By setting, for a generic continuous function $\phi:[0,+\infty)\rightarrow\RR$,
$$|\phi|_t = \max_{0\le s\le t}|\phi(s)|, \quad t\ge0,\qquad\qquad \|\phi\| = \sup_{t\ge0} |\phi(t)|\equiv\lim_{t\rightarrow\infty} |\phi|_t,$$ 
considering that
$$0\le s\le t \qquad \Rightarrow\qquad |\phi(s)|\le|\phi|_t,$$
denoting 
$$E_\aa^\lam = \sup_{t\ge0} |E_\aa(\lam t^\aa)|,$$
and taking into account (\ref{roa}) and (\ref{minor}), from (\ref{solf}) one derives:
\begin{eqnarray*}
|y(t)| &\le& |E_\aa(\lam t^\aa)||y_0| + \int_0^t (t-s)^{\aa-1}|E_{\aa\aa}(\lam(t-s)^\aa)||f(y(s))|\dd s\\
        &\le&  |E_\aa(\lam t^\aa)||y_0| + L\int_0^t (t-s)^{\aa-1}|E_{\aa\aa}(\lam(t-s)^\aa)||y(s)|\dd s\\
        &\le&  E_\aa^\lam \,|y_0| + L\rho_\aa^\lam |y|_t \le E_\aa^\lam \,|y_0| + \frac{1}2 |y|_t, \qquad \forall t\ge0.
\end{eqnarray*}
Consequently, 
$$\forall s\in[0,t]\,:\,|y(s)|\le E_\aa^\lam \,|y_0| + \frac{1}2 |y|_t\qquad \Rightarrow\qquad |y|_t\le 2E_\aa^\lam \,|y_0| .$$
Since $\lam\in\Lambda_\aa$, $E_\aa^\lam<\infty$. Therefore, by considering the limit for $t\rightarrow\infty$, it follows that
$$
\|y\| \le 2 E_\aa^\lam \,|y_0| <\infty.
$$
This clearly implies the stability of the origin. To prove its asymptotic stability, we observe, again from (\ref{solf}), that
\begin{eqnarray}\nonumber
\hat\ell &:=& \limsup_{t\rightarrow\infty} |y(t)| \le |y_0| \limsup_{t\rightarrow\infty}|E_\aa(\lam t^{\aa})| + L\limsup_{t\rightarrow\infty} \int_0^t (t-s)^{\aa-1}|E_{\aa\aa}(\lam(t-s)^\aa)||y(s)|\dd s\\ \label{ell}
&=&  L\limsup_{t\rightarrow\infty} \int_0^t (t-s)^{\aa-1}|E_{\aa\aa}(\lam(t-s)^\aa)||y(s)|\dd s,
\end{eqnarray}
due to the fact that (recall (\ref{talfa})) ~$\limsup_{t\rightarrow\infty}|E_\aa(\lam t^\aa)| =0.$~ 
Should we have $\hat\ell>0$, one would infer that 
$$\forall\eps>0 ~\exists\, T>0 ~ s.t.~ \forall s\ge T\,:\,|y(s)|\le\hat\ell+\eps,$$
i.e., $$\forall\eps>0 ~\exists\, T>0 ~ s.t.~ \sup_{s\ge T} |y(s)| \le \hat\ell+\eps.$$
Therefore,  since 
$$\limsup_{t\rightarrow\infty} \int_T^t(t-s)^{\aa-1}|E_{\aa\aa}(\lam(t-s)^\aa)||y(s)|\dd s ~\le~ (\hat\ell+\eps)\sup_{t\ge T} \int_T^t (t-s)^{\aa-1}|E_{\aa\aa}(\lam(t-s)^\aa)|\dd s,$$
from (\ref{ell}) it follows that:
$$\hat\ell \le L\left[ \|y\| \limsup_{t\rightarrow\infty} \int_0^T (t-s)^{\aa-1}|E_{\aa\aa}(\lam(t-s)^\aa)|\dd s
+(\hat\ell+\eps)\sup_{t\ge T} \int_T^t (t-s)^{\aa-1}|E_{\aa\aa}(\lam(t-s)^\aa)|\dd s\right].$$
Moreover, due to (\ref{talfa1}) we have that 
\begin{eqnarray*}
\limsup_{t\rightarrow\infty} \int_0^T (t-s)^{\aa-1}|E_{\aa\aa}(\lam(t-s)^\aa)|\dd s 
&=& \limsup_{t\rightarrow\infty} \int_{t-T}^t s^{\aa-1}|E_{\aa\aa}(\lam s^\aa)|\dd s \\
&=& T \limsup_{t\rightarrow\infty} t^{\aa-1}|E_{\aa\aa}(\lam t^\aa)|  ~=~ 0.
\end{eqnarray*}
 Consequently, by also taking into account that (see (\ref{roa})) 
 $$ \sup_{t\ge T} \int_T^t (t-s)^{\aa-1}|E_{\aa\aa}(\lam(t-s)^\aa)|\dd s\le \rho_\aa^\lam,$$
 one obtains the following inequality: 
 $$\hat\ell \le L\rho_\aa^\lam (\hat\ell+\eps) < \frac{\hat\ell+\eps}2.$$
Thus, by choosing any $\eps\le\hat\ell$, one has a contradiction. Therefore, we conclude that $\hat\ell=0$, which implies that
 $$\limsup_{t\rightarrow\infty} |y(t)|=0 \qquad \Longrightarrow\qquad \lim_{t\rightarrow\infty} y(t)=0,$$
 i.e., the origin is asymptotically stable.\,\QED

\subsection{The discrete case}\label{discrete}
Hereafter, we follow similar steps as those made above for discussing the continuous case. 
In particular, we shall hereafter consider the case of FHBVM$(k,s)$ methods with $k=s$ in (\ref{Butab})-(\ref{rankA}), i.e., the case of collocation methods. In fact, for such methods the Butcher matrix $A$ is nonsingular, according to (\ref{rankA}), and this property simplifies the arguments.\footnote{We shall consider the general case elsewhere.}
Consequently, at first let us consider the application of a FHBVM$(s,s)$ method for solving the linear test equation (\ref{test}) with timestep $t>0$. By setting
\begin{equation}\label{q}
q := \lam t^\aa,
\end{equation}
and $e=(1,\dots,1)^\top\in\RR^s$, standard arguments used for Runge-Kutta methods prove that 
\begin{equation}\label{RKtest}
Y = ey_0 +q AY, \qquad y_1 =  y_0 + \frac{q}{\Gamma(\aa+1)}\bfb^\top Y, 
\end{equation}
from which one readily derives that,\footnote{Hereafter, for a FHBVM$(s,s)$ method, $I$ will denote the identity matrix of dimension $s$, for the sake of brevity.}  assuming  matrix $I-qA$ be nonsingular,
\begin{equation}\label{Rq}
y_1 = \left( 1 + \frac{q}{\Gamma(\aa+1)}\bfb^\top(I-q A)^{-1}e\right) y_0 =: R_s^\aa(q) y_0.
\end{equation}
A straightforward comparison of (\ref{sol}) and (\ref{Rq}) then shows that
\begin{equation}\label{approE}
R_s^\aa(q)\approx E_\aa(q).
\end{equation}
As is usual, we shall define the stability region of the discrete method as:
\begin{equation}\label{D}
\D_s^\aa := \left\{ q\in\CC\,:\, |R_s^\aa(q)|<1\right\}.
\end{equation}
Similarly as in \cite{WZ2023}, we give the following definition.
\begin{defi}\label{Astab}
The FHBVM$(s,s)$ method is said to be $A$-stable when, with reference to (\ref{Ea0}), 
$$\forall z\in\CC \,:\,|E_\aa(z)|<1 ~\Rightarrow~ |R_s^\aa(z)|<1, \qquad \Leftrightarrow\qquad
\Lambda_\aa\subseteq\D_s^\aa.$$
\end{defi}
\begin{rem} 
Taking into account (\ref{q}), in order for matrix $(I-qA)\equiv (I-\lam t^\aa A)$ to be nonsingular for all $\lam\in\Lambda_\aa$ and $t>0$, matrix $A$ must satisfy:\footnote{As is usual, $\sigma(A)$ denotes the spectrum of matrix $A$.}
\begin{equation}\label{muA}
\mu\in\sigma(A) \qquad \Rightarrow \qquad |\arg(\mu)|<\aa\frac{\pi}2.
\end{equation}
We have numerically verified that (\ref{muA}) holds true for 
$$s=1,\dots,50, \qquad\mbox{and}\qquad  \aa=0.001,0.002,\dots,0.998,0.999,$$ 
so that, hereafter, we shall assume (\ref{muA}) holds true in general. This fact is consistent with the limit case when $\aa\rightarrow1$, where  FHBVM$(s,s)$ becomes the $s$-stage Gauss-Legendre collocation method.
\end{rem}

\medskip 

As in the continuous case, the results concerning the test equation (\ref{test}) can be extended to the nonlinear problem (\ref{testf}), with $f(0)=0$ and $f$ Lipschitz, provided that the Lipschitz constant $L$ is sufficiently small.\footnote{As in the continuous case, alternative, though essentially equivalent, requirements on $f$ can be also considered.} In fact, now the equations in (\ref{RKtest}) take the form
\begin{equation}\label{RKtestf}
Y = ey_0+t^\aa A(\lam Y+f(Y)), \qquad y_1 = y_0 + \frac{t^\aa}{\Gamma(\aa+1)}\bfb^\top(\lam Y+f(Y)).
\end{equation}
Some algebra allows us to derive from (\ref{RKtestf}), by considering (\ref{q}), (\ref{Rq}), and assuming that matrix $(I-qA)$ is nonsingular:
\begin{eqnarray}\label{Yq}
Y     &=& (I-qA)^{-1}ey_0 + t^\aa(I-qA)^{-1}Af(Y),\\[1mm] \label{y1q}
y_1 &=& R_s^\aa(q)y_0 + \frac{t^\aa}{\Gamma(\aa+1)}\bfb^\top (I-qA)^{-1} f(Y).
\end{eqnarray}

The following result then holds true.

\begin{theo}\label{serveRq} Assume that $\lam\in\Lambda_\aa$ and, moreover, for any suitable norm $\|\cdot\|$:
\begin{equation}\label{L0}
L \max(\|\bfb\|,\|A\|)\sup_{t\ge0}t^\aa\|(I-q A)^{-1}\|\equiv L \max(\|\bfb\|,\|A\|)\sup_{t\ge0}t^\aa\|(I-\lam t^\aa A)^{-1}\|<\frac{1}2.
\end{equation}
Then, (\ref{Rq}) and (\ref{y1q}) have the same asymptotic behavior, as $t\rightarrow\infty$.
\end{theo}
\proof
To begin with, from (\ref{Yq}) one has
$$\|Y\|\le \|(I-qA)^{-1}e\||y_0| + L t^\aa\|(I-qA)^{-1}A\|\|Y\|$$
and, therefore, by virtue of (\ref{L0}):
$$\|Y\|\le 2\|(I-qA)^{-1}e\||y_0|.$$
Consequently, from (\ref{y1q}), again by using (\ref{L0}), one derives that: 
$$\left| y_1 - R_s^\aa(q)y_0\right| \le \frac{1}{\Gamma(\aa+1)}\|(I-qA)^{-1}e\||y_0|.$$
Considering that, by virtue of (\ref{q}),  $\|(I-qA)^{-1}e\|=O(t^{-\aa})$, as $t\rightarrow\infty$,  the statement then follows.\,\QED
\bigskip

From the result of Theorem~\ref{serveRq}, one infers that:
\begin{itemize}
\item the better $R_s^\aa(q)$ approximates $E_\aa(q)$, the better the numerical approximation is (see (\ref{approE}));

\item it is known that, for $\lam\in\Lambda_\aa$,  $E_\aa(\lam t^\aa)=O(t^{-\aa})$ as $t\rightarrow\infty$. On the other hand, from (\ref{Rq}) one derives that 
\begin{equation}\label{Rinf}
R_s^\aa(\infty)~:=~\lim_{t\rightarrow\infty}R_s^\aa(\lam t^\aa) ~=~ 1-\frac{\bfb^\top A^{-1}e}{\Gamma(\aa+1)}.
\end{equation}
\end{itemize}

\begin{figure}[p]
\centering
\includegraphics[width=8cm]{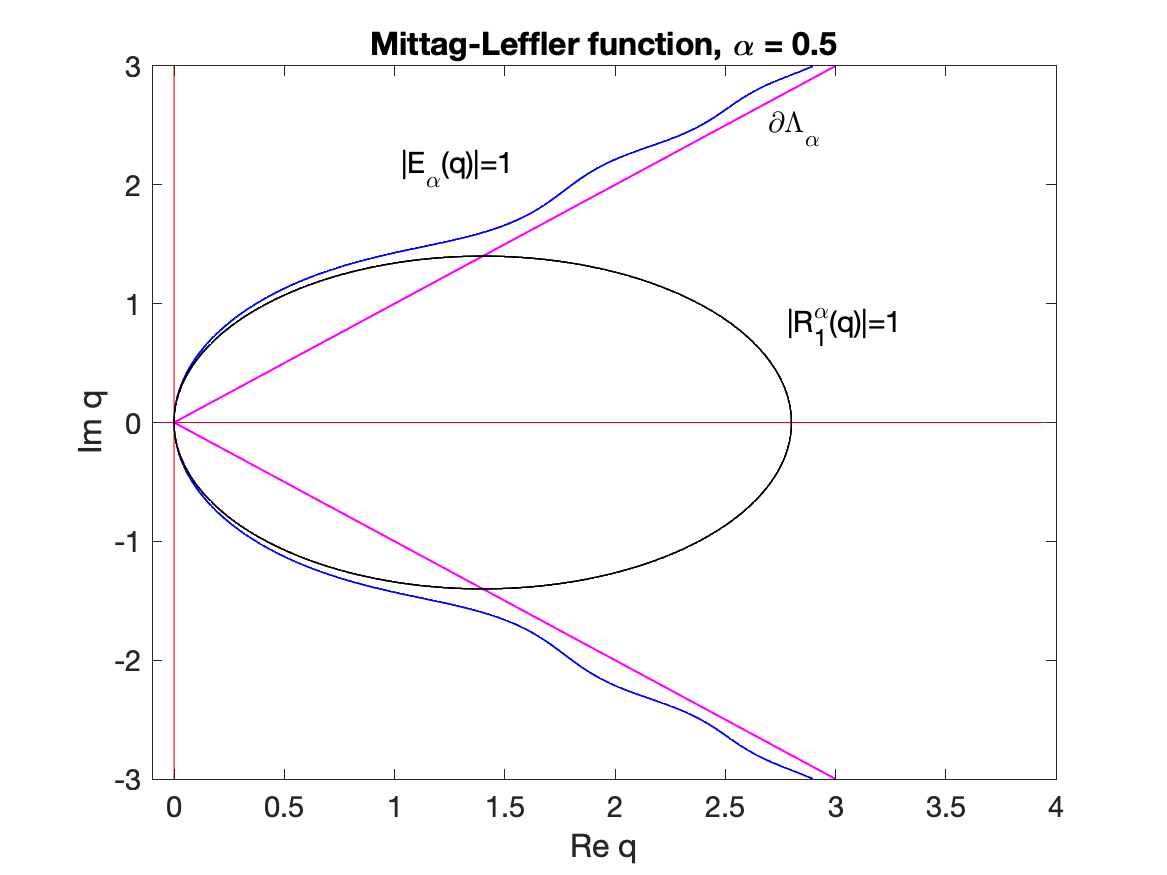}

\includegraphics[width=8cm]{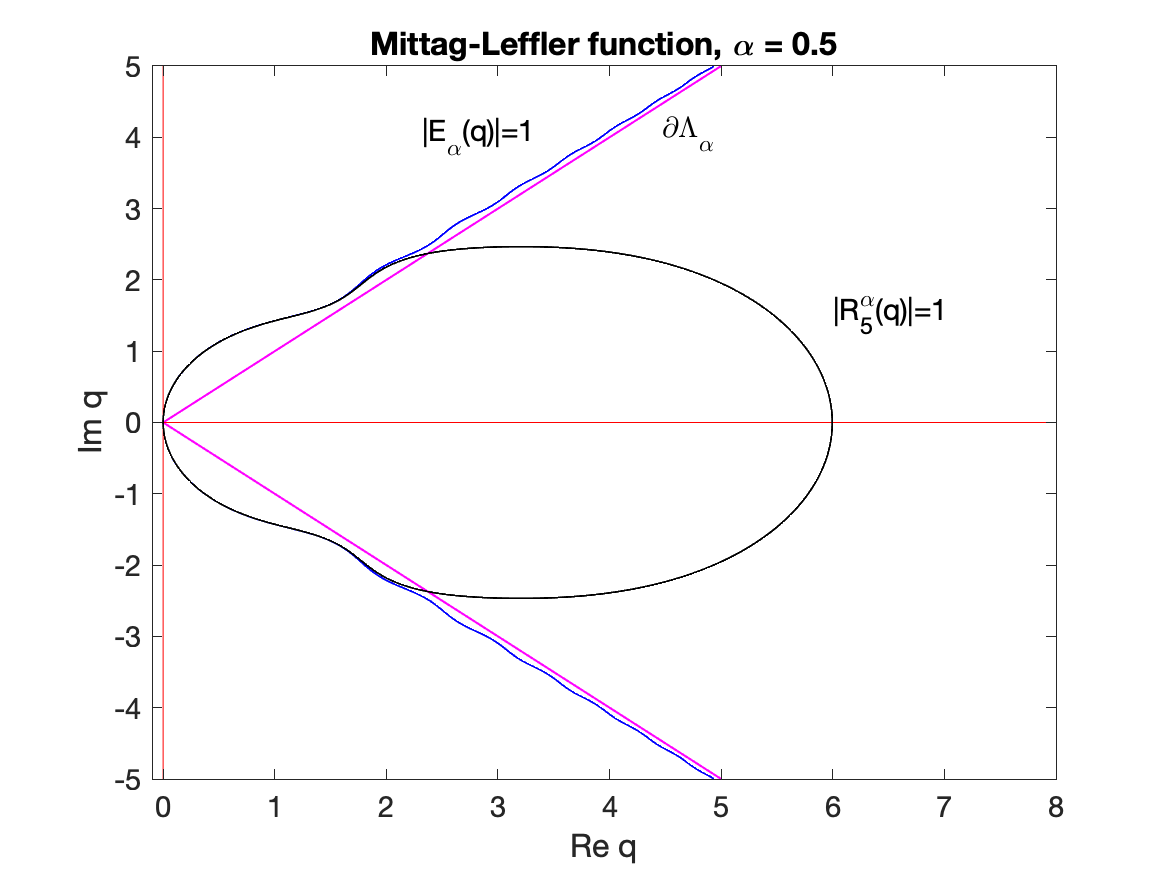}

\includegraphics[width=8cm]{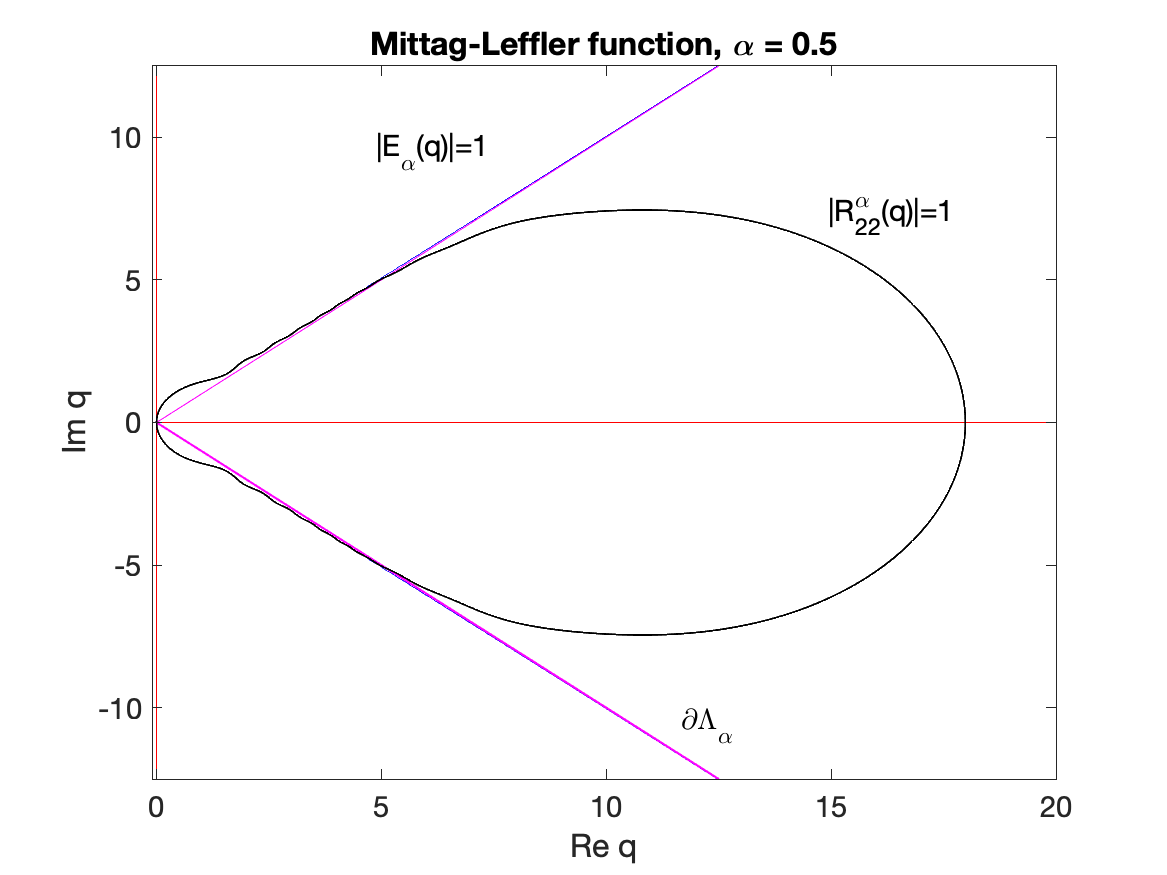}
\caption{Boundary of the stability region of FHBVM$(s,s)$, $s=1$ (upper plot), $s=5$ (middle plot), $s=22$ (lower plot), $\aa=0.5$.}
\label{seq1to22}
\end{figure}

\smallskip
Concerning the first point, in the three plots of Figure~\ref{seq1to22} we report, for the case $\aa=0.5$:\footnote{For computing the Mittag-Leffler function, we have used the function {\tt ml.m} in \cite{MLfun}.}
\begin{enumerate}
\item the level set $|E_\aa(q)|=1$;
\item $\partial\Lambda_\aa$, i.e., the boundary of $\Lambda_\aa$; 
\item the level set $|R_s^\aa(q)|=1$, $s=1,5,22$, which is the boundary of the corresponding $\D_s^\aa$ region (\ref{D}) (i.e., the unbounded outer one).
\end{enumerate}
From these plots, one deduces that all methods are $A$-stable and, the larger the parameter $s$ is, the better $R_s^\aa(q)$ approximates $E_\aa(q)$, as is expected.

\smallskip

\begin{figure}[t]
\centering
\includegraphics[width=10cm]{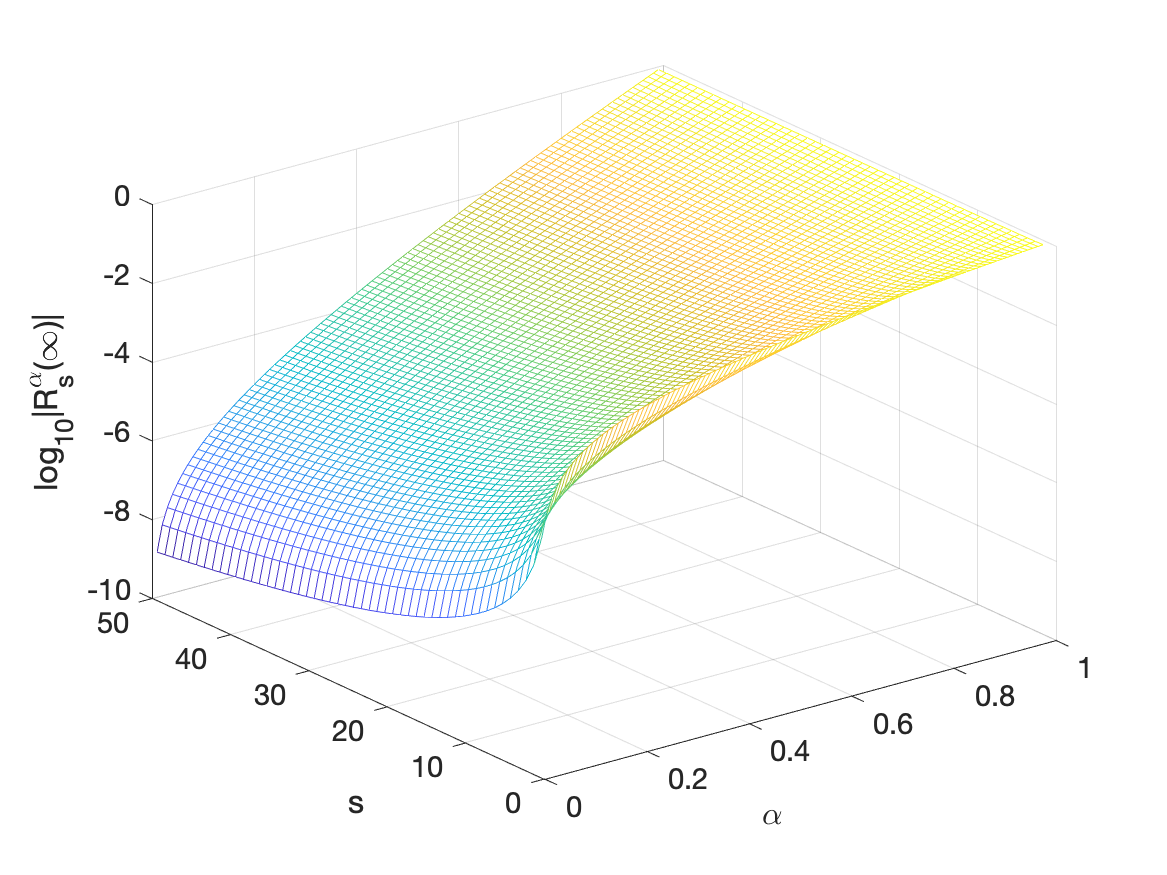}
\caption{Absolute value of the parameter (\ref{Rinf}), for $\aa\in(0,1)$ and $s=1,\dots,50$.}
\label{rinfig}
\end{figure}

Concerning (\ref{Rinf}), in Figure~\ref{rinfig} we plot $\log_{10}|R_s^\aa(\infty)|$ for $\aa\in(0,1)$ and $s=1,\dots,50$. From this figure, one deduces that $$|R_s^\aa(\infty)|\rightarrow 1, \quad \mbox{as} \quad \aa\rightarrow 1,$$ which is expected, since in such a case the method tends to the corresponding $s$-stage Gauss-Legendre collocation formula. Similarly,  as $\aa\rightarrow0$, $R_s^\aa(\infty)$ tends to 0. Moreover, $|R_s^\aa(\infty)|$ is a decreasing function of $s$: i.e., the larger the parameter $s$ is, the more consistent the behavior of the method is with the continuous case. This fact is in accord with the expected spectral accuracy of the methods, for larger values of $s$.

\section{The mixed stepsize selection}\label{appros}

We now provide full details for the mixed stepsize selection sketched in the introduction. 
To begin with, with reference to (\ref{h}), (\ref{hi}), and (\ref{rnu}), let us define the mesh:
\begin{equation}\label{mesh}
\hat t_0 = 0, \quad \hat t_i = \hat t_{i-1}+h_i, \qquad i=1,\dots,\nu, \qquad t_j = j h \equiv j\frac{T}N, \qquad j=n,\dots,N.
\end{equation}
Clearly, (\ref{mesh}) reduces to a pure graded mesh, when $\nu>n=N$, or to a purely uniform mesh, when $\nu=n$, so that the approach encompasses, as particular cases, the previous ones. Further, let us denote, for $c\in[0,1]$,
\begin{equation}\label{restri}
\hat y_i(ch_i) \equiv y(\hat t_{i-1}+c h_i), \quad i=1,\dots,\nu, \qquad y_j(ch)\equiv y(t_{j-1}+ch), \quad j=n+1,\dots,N,
\end{equation}
the restrictions of the solution of (\ref{fde}) to the respective sub-intervals $[\hat t_{i-1},\hat t_i]$ and $[t_{j-1},t_j]$. Consequently, by using similar steps as in \cite{BGI2024}, for $i=1,\dots,\nu$\,  one obtains, by setting (see (\ref{solfde})) $\hat T_\ell^i(ch_i)\equiv T_\ell(\hat t_{i-1}+ch_i)$:
\begin{eqnarray}\nonumber
\lefteqn{\hat y_i(ch_i)~=~ y(\hat t_{i-1}+ch_i)~=~T_\ell(\hat t_{i-1}+ch_i) ~+~ I^\aa f(\hat y(\hat t_{i-1}+ch_i)) }\\ \nonumber
&=& \hat T_\ell^i(ch_i) ~+~ \frac{1}{\Gamma(\aa)} \left[ \sum_{\mu=1}^{i-1} \int_{\hat t_{\mu-1}}^{\hat t_\mu} (\hat t_{i-1}+ch_i-x)^{\aa-1}f(y(x))\dd x \right.\\ \nonumber &&\left.\qquad\quad~\, +~  \int_{\hat t_{i-1}}^{\hat t_{i-1}+ch_i} (\hat t_{i-1}+ch_i-x)^{\aa-1}f(y(x))\dd x \right]\\
\nonumber
&=& \hat T_\ell^i(ch_i) ~+~ \frac{1}{\Gamma(\aa)} \left[ \sum_{\mu=1}^{i-1} 
h_\mu^\aa\int_0^1\left(\frac{r^{i-\mu}-1}{r-1}+cr^{i-\mu}-\tau\right)^{\aa-1}f(\hat y_\mu(\tau h_\mu))\dd \tau  \right.\\ 
&&\qquad\quad~\, \left.  +~h_i^\aa\int_{0}^c (c-\tau )^{\aa-1}f(\hat y_i(\tau h_i))\dd \tau \right], \qquad c\in[0,1].\label{hyi}
\end{eqnarray}
Similarly, for $j=n+1,\dots,N$, by setting $T_\ell^j(ch) \equiv T_\ell( t_{j-1}+ch)$, one derives:
\begin{eqnarray}\nonumber
\lefteqn{y_j(ch)~=~ y(t_{j-1}+ch)~=~T_\ell( t_{j-1}+ch) ~+~ I^\aa f(y(t_{j-1}+ch))}\\ \nonumber
&=& T_\ell^j(ch) ~+~ \frac{1}{\Gamma(\aa)} \left[ \sum_{i=1}^\nu \int_{\hat t_{i-1}}^{\hat t_i} (t_{j-1}+ch-x)^{\aa-1}f(y(x))\dd x \right.\\ \nonumber
&&\left. +~ \sum_{\mu=n+1}^{j-1} \int_{t_{\mu-1}}^{t_\mu} (t_{j-1}+ch-x)^{\aa-1}f(y(x))\dd x ~+~
\int_{t_{j-1}}^{t_{j-1}+ch} (t_{j-1}+ch-x)^{\aa-1}f(y(x))\dd x\right]
\\ \nonumber
&=& T_\ell^j(ch) ~+~ \frac{1}{\Gamma(\aa)} \left[ \sum_{i=1}^\nu 
h_i^\aa\int_0^1\left( \frac{(j-1+c)(r^\nu-1)-n(r^{i-1}-1)}{nr^{i-1}(r-1)}-\tau\right)^{\aa-1}f(\hat y_i(\tau h_i))\dd \tau  \right.\\ \nonumber
&&\left.  +~ h^\aa\sum_{\mu=n+1}^{j-1}     
\int_0^1\left(j-\mu+c-\tau\right)^{\aa-1}f(y_\mu(\tau h))\dd \tau + h^\aa\int_{0}^c (c-\tau )^{\aa-1}f(y_j(\tau h))\dd \tau \right], \\ 
&&\qquad c\in[0,1].\quad \label{yj}
\end{eqnarray}

Following steps similar to those used in Section~\ref{fhbvm_sec}, a quasi-polynomial approximation can then be derived by considering that (\ref{fde}) can be reformulated as the following system of local FDEs (see (\ref{restri})):
\begin{eqnarray}\nonumber
\hat y_i^{(\aa)}(ch_i) &=& f(\hat y_i(ch_i)), \qquad i=1,\dots,\nu, \qquad \hat y_1^{(\mu)}(0) = y_0^\mu, \qquad \mu=0,\dots,\ell-1,\\[2mm]
y_j^{(\aa)}(ch) &=& f(y_j(ch)), \qquad j=n+1,\dots,N,\qquad c\in[0,1].\label{sys1}
\end{eqnarray}
By considering the expansions of the local vector fields in (\ref{sys1}) along the orthonormal Jacobi basis (\ref{orto}),
\begin{equation}\label{exp1}
f(\hat y_i(ch_i)) = \sum_{\iota\ge 0} P_\iota(c)\gamma_\iota^i(\hat y_i), \qquad 
f(y_j(ch)) = \sum_{\iota\ge 0} P_\iota(c)\gamma_\iota(y_j), \qquad c\in[0,1],
\end{equation}
where, for any suitable given function $z$,
\begin{equation}\label{gammai}
\gamma_\iota^i(z) = \int_0^1\omega(c)P_\iota(c)f(z(c h_i))\dd c, \qquad \gamma_\iota(z) = \int_0^1\omega(c)P_\iota(c)f(z(c h))\dd c, \qquad \iota=0,1,\dots,
\end{equation}
an approximate polynomial vector field can be obtained by truncating the infinite series in (\ref{exp1}) to finite sums with $s$ terms. In so doing, one derives
local approximations $$\hat \sigma_i(ch_i)\approx \hat y_i(ch_i), \quad i=1,\dots,\nu, \qquad  \sigma_j(ch)\approx y_j(ch), \quad j=n+1,\dots,N, \qquad c\in[0,1],$$
formally satisfying the approximate system of FDEs:
\begin{eqnarray}\nonumber
\hat \sigma_i^{(\aa)}(ch_i) &=& \sum_{\iota = 0}^{s-1} P_\iota(c)\gamma_\iota^i(\hat \sigma_i), \qquad i=1,\dots,\nu, \qquad \hat \sigma_1^{(\mu)}(0) = y_0^\mu, \qquad \mu=0,\dots,\ell-1,\\
\sigma_j^{(\aa)}(ch) &=&  \sum_{\iota = 0}^{s-1} P_\iota(c)\gamma_\iota(\sigma_j), \qquad j=n+1,\dots,N,\qquad c\in[0,1],\label{sys2}
\end{eqnarray}
with the Fourier coefficients defined according to (\ref{gammai}). A corresponding expression of the solution is then derived from (\ref{hyi}) and (\ref{yj}), respectively, by formally replacing the vector fields in (\ref{sys1}) with the corresponding ones in (\ref{sys2}). In so doing, for $i=1,\dots,\nu$, one obtains:
\begin{eqnarray*}
\hat \sigma_i(ch_i)&=& \hat T_\ell^i(ch_i) ~+~ \frac{1}{\Gamma(\aa)} \left[ \sum_{\mu=1}^{i-1} 
h_\mu^\aa\int_0^1\left(\frac{r^{i-\mu}-1}{r-1}+cr^{i-\mu}-\tau\right)^{\aa-1}\sum_{\iota=0}^{s-1} P_\iota(\tau) \gamma_\iota^\mu(\hat \sigma_\mu)\dd \tau \right.\\  
&& \left.  +~ h_i^\aa\int_{0}^c (c-\tau )^{\aa-1}\sum_{\iota=0}^{s-1} P_\iota(\tau) \gamma_\iota^i(\hat \sigma_i)\dd \tau \right]\\
&\equiv& \hat\phi_i(c) ~+~ h_i^\aa\sum_{\iota=0}^{s-1} I^\aa P_\iota(c)\gamma_\iota^i(\hat\sigma_i),\qquad c\in[0,1], 
\end{eqnarray*}
with the {\em memory term} $\hat\phi_i(c)$ defined as 
$$
\hat\phi_i(c) ~=~ \hat T_\ell^i(ch_i) ~+~ \sum_{\mu=1}^{i-1}h_\mu^\aa \sum_{\iota=0}^{s-1} J_\iota^\aa\left(\frac{r^{i-\mu}-1}{r-1}+cr^{i-\mu}\right)\gamma_\iota^\mu(\hat\sigma_\mu),
$$
having set
\begin{equation}\label{Jint}
J_\iota^\aa(x) ~=~ \frac{1}{\Gamma(\aa)}\int_0^1 (x-\tau)^{\aa-1}P_\iota(\tau)\dd\tau,
\end{equation}
and with $I^\aa P_\iota(c)$ the Riemann-Liouville integral of $P_\iota(c)$ (see (\ref{solfde})). Similarly, from (\ref{yj}), for $j=n+1,\dots,N$, one derives:
$$
\sigma_j(ch) ~=~ \phi_j(c) ~+~ h^\aa\sum_{\iota=0}^{s-1} I^\aa P_\iota(c) \gamma_\iota(\sigma_j), \qquad c\in[0,1],
$$
with the memory term now given by (see (\ref{Jint})):
\begin{eqnarray*}
\phi_j(c) &=& T_\ell^j(ch) ~+~ \sum_{i=1}^\nu h_i^\aa \sum_{\iota=0}^{s-1} J_\iota^\aa\left( \frac{(j-1+c)(r^\nu-1)-n(r^{i-1}-1)}{nr^{i-1}(r-1)}\right)\gamma_\iota^i(\hat\sigma_i)\\ 
&&~+~h^\aa\sum_{\mu=n+1}^{j-1}\sum_{\iota=0}^{s-1} J_\iota^\aa(j-\mu+c)\gamma_\iota(\sigma_\mu).
\end{eqnarray*}

\section{Dicretization}\label{discre}

Clearly (see also \cite{BGI2024}), the discrete solution is completely known once the Fourier coefficients in (\ref{sys2}) are computed: in particular those at the current step, since the previous ones have already been obtained. This amounts to solving discrete problems either in the form
\begin{equation}\label{hgi}
\gamma_\iota^i(\hat\sigma_i) = \int_0^1\omega(c) P_\iota(c) f\left( \hat\phi_i(c)+h_i^\aa\sum_{\mu=0}^{s-1} I^\aa P_\mu(c)\gamma_\mu^i(\hat\sigma_i)\right)\dd c, \qquad \iota = 0,\dots,s-1,
\end{equation}
when we approximate the solution on the $i$th sub-interval of the graded mesh, or
\begin{equation}\label{gj}
\gamma_\iota(\sigma_j) = \int_0^1\omega(c) P_\iota(c) f\left( \phi_j(c)+h^\aa\sum_{\mu=0}^{s-1} I^\aa P_\mu(c)\gamma_\mu(\sigma_j)\right)\dd c, \qquad \iota = 0,\dots,s-1,\end{equation}
when we approximate the solution on the $j$th sub-interval of the uniform mesh.

As recalled in Section~\ref{fhbvm_sec},  using the Gauss-Jacobi quadrature formula of order $2k$ for approximating the integrals in (\ref{hgi})-(\ref{gj}) results into a FHBVM$(k,s)$ method. This has the advantage of leaving unaltered the (block) dimension $s$ of the discrete problems, whichever is the chosen value of $k$, and, moreover, only requires to compute the involved functions at the quadrature abscissae $c_\rho$, $\rho=1,\dots,k$. This, in turn, implies that we only need to evaluate the following integrals, for 
all ~$\iota = 0,\dots,s-1$,~ and ~$\rho=1,\dots,k$:\footnote{We refer to \cite{BGI2024} and \cite{ABI2022} for the numerical evaluation of such integrals.}
\begin{eqnarray} \label{I1}
I^\aa P_\iota(c_\rho),&& \\ \label{I2}
J_\iota^\aa(j+c_\rho),&&  \qquad j=1,\dots,N-n-1, \\ \label{I3}
J_\iota^\aa\left(\frac{r^i-1}{r-1}+c_\rho r^i\right),&& \qquad i=1,\dots,\nu-1, \\ \label{I4}
J_\iota^\aa\left(\frac{(j+c_\rho)(r^\nu-1)-n(r^i-1)}{nr^i(r-1)}\right),&&\qquad i=0,\dots,\nu-1,\quad j=n,\dots,N-1.
\end{eqnarray}
The code {\tt fhbvm} \cite{BGI2024}, available at the URL \cite{fhbvm}, basically computes the integrals (\ref{I1}), and either (\ref{I2}) or (\ref{I3}), depending on the mesh used (uniform or graded, respectively). On the other hand, the new mixed-mesh strategy here described requires to compute both  the integrals in (\ref{I2}) and (\ref{I3}), along with those in (\ref{I4}): these latter integrals, in turn, are computed following the approximation procedure explained in \cite{BGI2024} for computing the integrals (\ref{I2}) and (\ref{I3}). Such modifications have been implemented in the Matlab\cpr code {\tt fhbvm2},  based on the FHBVM(22,22) method (i.e., the collocation method corresponding to $k=s=22$). The code, made available at the same URL \cite{fhbvm} of the code {\tt fhbvm}, will be used for the numerical tests in the next section.

\section{Numerical tests}\label{numtest}

Here, we report some comparisons among the following Matlab\cpr codes:

\begin{itemize}
\item {\tt flmm2} \cite{Garr18}, by selecting the BDF2 method;\footnote{We have used the
Revision: 1.0 - Date: June 27 2014, of the code {\tt flmm2}.}
\item {\tt flmm2} \cite{Garr18}, by selecting the trapezoidal rule;
\item {\tt fhbvm} \cite{BGI2024}; 
\item the code {\tt fhbvm2} here described.\footnote{Likewise the code {\tt fhbvm}, also {\tt fhbvm2} uses the Matlab\cpr suite OPQ, containing companion codes of \cite{Gautschi2004}.} 
\end{itemize}
For the {\tt flmm2} code the  parameters {\tt tol=1e-15} and {\tt itmax=1000} have been used.
All numerical tests have been done on a  8-core M2-Silicon based computer with 16GB of shared memory, using Matlab\cpr  R2024b. The accuracy of the numerical solutions is measured in terms of  {\em mixed error significant computed digits} ({\tt mescd}), defined as \cite{testset}: 
$${\tt mescd} := \max\left\{0,-\log_{10} \max_i \|(\bar y_i - y_i)./(1+|\bar y_i|)\|_\infty\right\},$$ 
where $y_i$ is the computed approximation at the $i$th mesh point, $\bar y_i$ is the corresponding reference solution, $|\cdot|$ is the vector with the absolute values, and $./$ is the component-wise division.

Comparisons among the methods are done through a corresponding {\em Work-Precision Diagram (WPD)} \cite{testset}, where the execution time (in sec) is plotted against accuracy. The elapsed time has been measured through the {\tt tic} and {\tt toc} functions of Matlab\cpr.\footnote{This is a difference w.r.t. \cite{BBBI2024,BGI2024,BGI2025}, where the function {\tt cputime} was used, instead.}

The calling sequence of the new code {\tt fhbvm2} is

\medskip
\centerline{\tt [t,y,etim] = fhbvm2( fun, y0, T, N, n, nu )}

\smallskip
\no where, in input:

\begin{itemize}
\item {\tt fun}\, contains the identifier of the function implementing the vector field (also in vector mode) and its Jacobian. Moreover, when called without arguments, the function must return the order $\aa$ of the fractional derivative;

\item {\tt y0}\, is an $\ell\times m$ matrix containing the initial conditions in (\ref{fde}), in the given order;

\item {\tt T}\, is the final integration time $T$;

\item {\tt N, n, nu}\, are the parameters $N,n,\nu$ for the mixed mesh (\ref{h})--(\ref{rnu}) seen above, whereas the parameter $r$ of the graded mesh is automatically set as follows:
$$r=\left\{ \begin{array}{ccl} 2, &\mbox{if}& n=1,\\[3mm]
\displaystyle\frac{n}{n-1}, &\mbox{if}& n>1.\end{array}\right.$$
This choice of the parameter, in turn, guarantees that the last timestep in the graded mesh satisfies: 
$$h_\nu \equiv h_1 r^{\nu-1} = \left\{
\begin{array}{ccc}\displaystyle
\frac{h}{r(1-r^{-\nu})},&\mbox{if} &n=1,\\[4mm]
 \displaystyle \frac{h}{1-r^{-\nu}},&\mbox{if} &n>1,\end{array}\right.$$ 
with $h$ the timestep used in the uniform mesh. 
One easily verifies that $h_\nu<h$, when $n=1$, whereas $h_\nu$ is sligthly larger than $h$, when $n>1$: in this latter case, if necessary, the input value of $\nu$ is increased, in order to have $h_\nu\le 1.1\cdot h$.
For the sake of completeness, we recall that the initial timestep of the graded mesh, $h_1$,  is computed from (\ref{rnu}):  a larger value of $\nu$ is thus required to reduce it. 
In output, the code provides:
\begin{itemize}

\item {\tt t, y}\, which contain the mesh and the computed solution, respectively;

\item {\tt etim}\, containing the execution time.

\end{itemize}
\end{itemize}

\subsection{Problem~1}\label{p1}
\begin{figure}[ph]
\centering
\includegraphics[width=8cm]{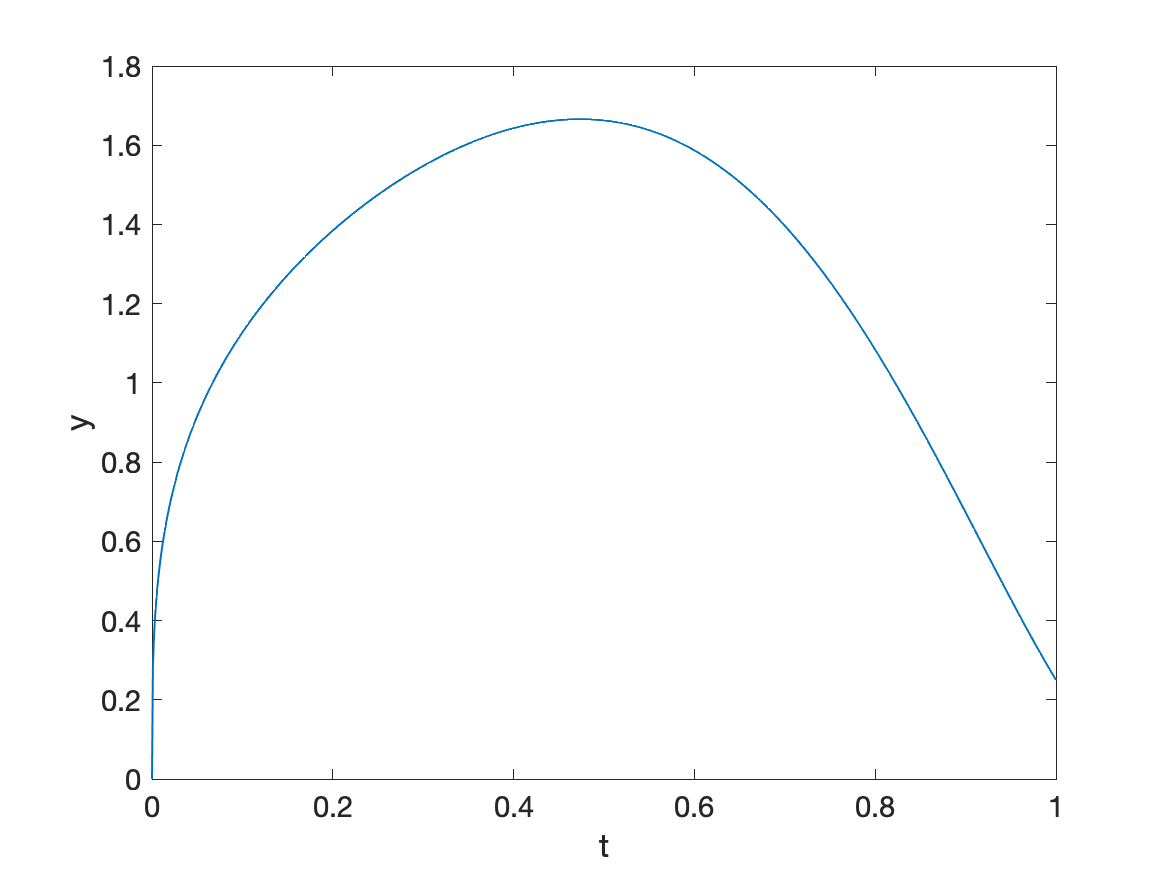}~\includegraphics[width=8cm]{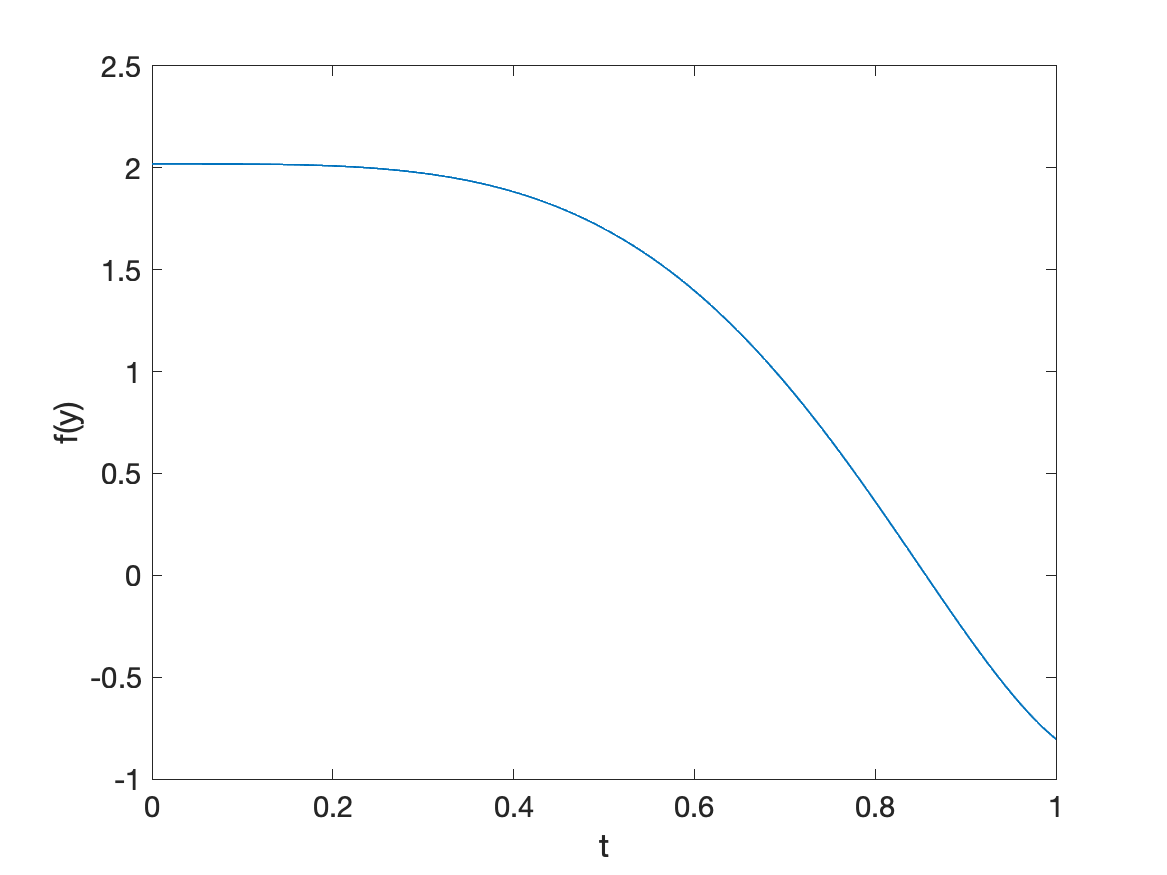}
\caption{Problem (\ref{ex1}), solution (left-plot) and vector field (right-plot).}
\label{ex1fig}

\bigskip\bigskip
\includegraphics[width=10cm]{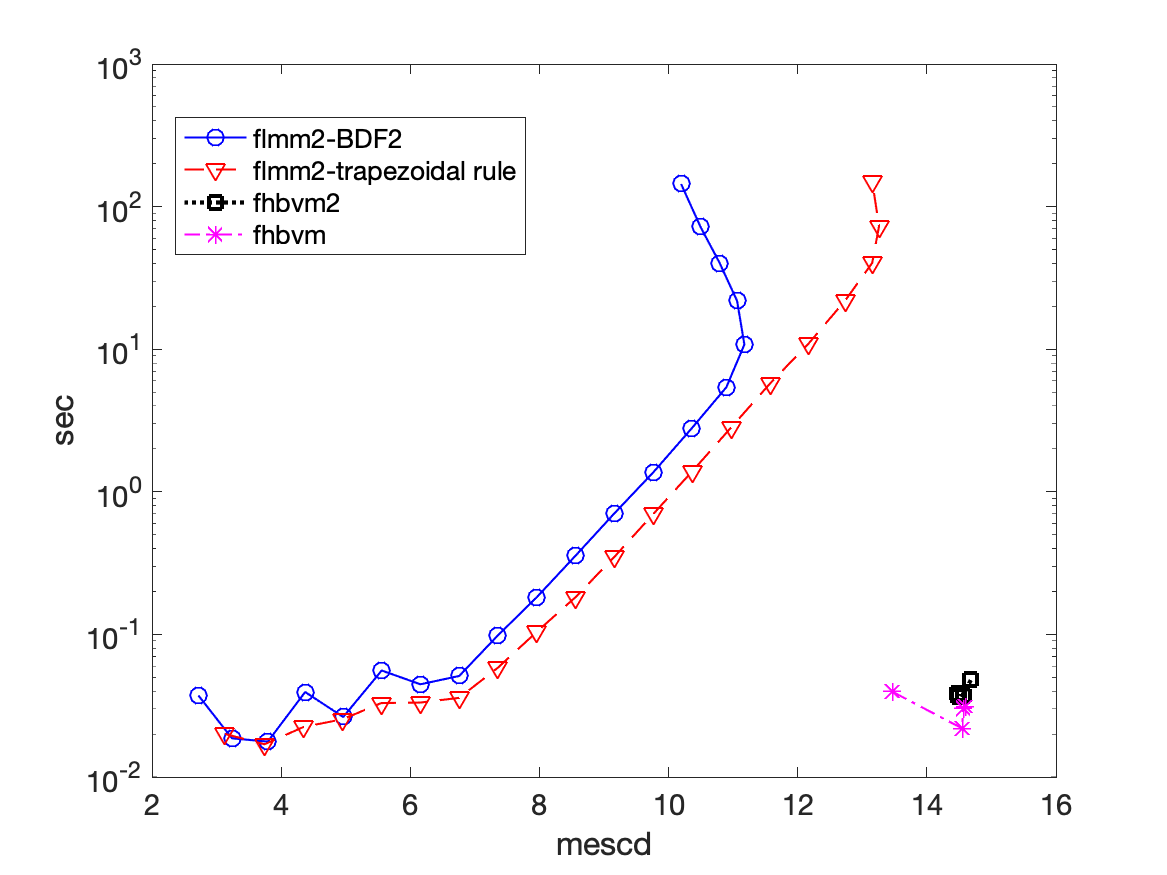}
\caption{Problem (\ref{ex1}), WPD.}
\label{ex1wpd}
\end{figure}

The first test problem  \cite{Garr18} is:
\begin{eqnarray}\nonumber
y^{(\aa)}(t) &=& -y^{3/2}(t) + \frac{8!}{\Gamma(9-\aa)}t^{8-\aa}-3\frac{\Gamma(5+\aa/2)}{\Gamma(5-\aa/2)}t^{4-\aa/2} +\left(\frac{3}2t^{\aa/2}-t^4\right)^3+\frac{9}4\Gamma(\aa+1),\\[1mm]
&& \qquad t\in[0,1], \qquad y(0)=0, \label{ex1}
\end{eqnarray}
whose solution is given by
\begin{equation}\label{ex1_sol}
y(t) = t^8-3t^{4+\aa/2}+\frac{9}4t^\aa.
\end{equation}
As in \cite{BGI2024}, we consider the value $\aa=0.3$. Even though the solution (\ref{ex1_sol}) is not smooth at the origin (see the left-plot in Figure~\ref{ex1fig}), the vector field turns out to be very smooth (as is shown in the right-plot in the same figure). Consequently, for this problem a uniform mesh is appropriate. For building the WPD relative to this problem, we use the codes with the following parameters:
\begin{itemize}
\item {\tt flmm2} (both using the trapezoidal rule and the BDF2 method): $h=0.1\cdot 2^{-i}$, $i=1,\dots,20$;
\item {\tt fhbvm}: $M=2,3,4,5$;\footnote{We recall that, for the code {\tt fhbvm}, $T/M$ is the used timestep, if a uniform mesh is selected, or is approximately equal to the last timestep, if a graded mesh is used. The choice between the two types of mesh is automatically done by the code \cite{BGI2024}.}
\item {\tt fhbvm2}: $\nu=n=1$, $N=2,3,4,5$.
\end{itemize}
Since a uniform mesh is used, we expect a similar performance for {\tt fhbvm} and {\tt fhbvm2}. This is indeed confirmed by the corresponding WPD, which is shown in Figure~\ref{ex1wpd}. From this figure, in fact, we can conclude that:
\begin{itemize}
\item {\tt flmm2} using BDF2 reaches about 11 significant digits of precision in about 11 sec. Further stepsize reductions do not improve the accuracy but only increase the execution time;

\item {\tt flmm2} using the trapezoidal rule reaches about 13 significant digits of precision in about 73 sec. As in the previous case, further stepsize reductions do not improve accuracy but only increase the execution time;

\item {\tt fhbvm} and {\tt fhbvm2} both provide a uniform accuracy of about 15 digits, independently of the stepsize used, with a negligible execution time (less than  $10^{-1}$ sec). 
\end{itemize}

\subsection{Problem~2}\label{p2}
\begin{figure}[ph]
\centering
\includegraphics[width=12cm]{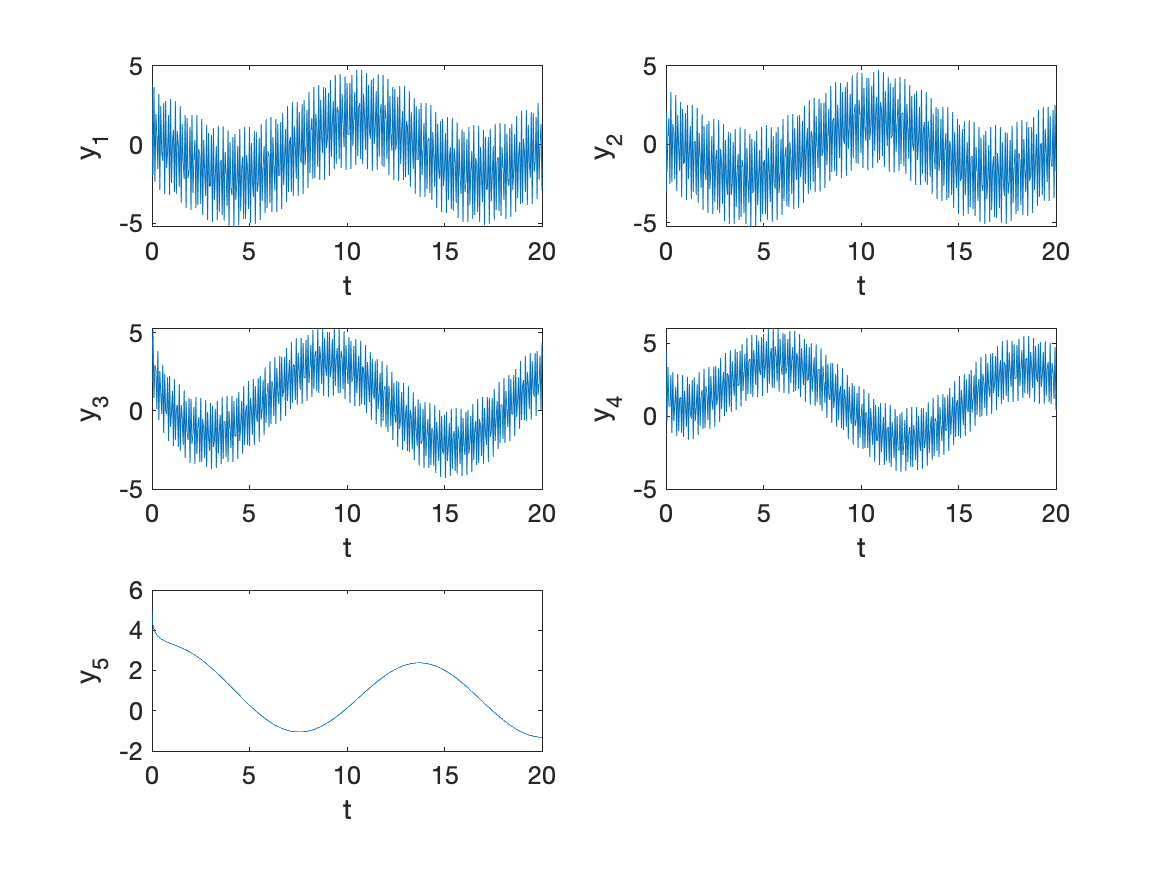}
\caption{Problem (\ref{ex3})-(\ref{ex3_Ay0}), components of the solution.}
\label{ex3fig}

\bigskip\bigskip
\includegraphics[width=10cm]{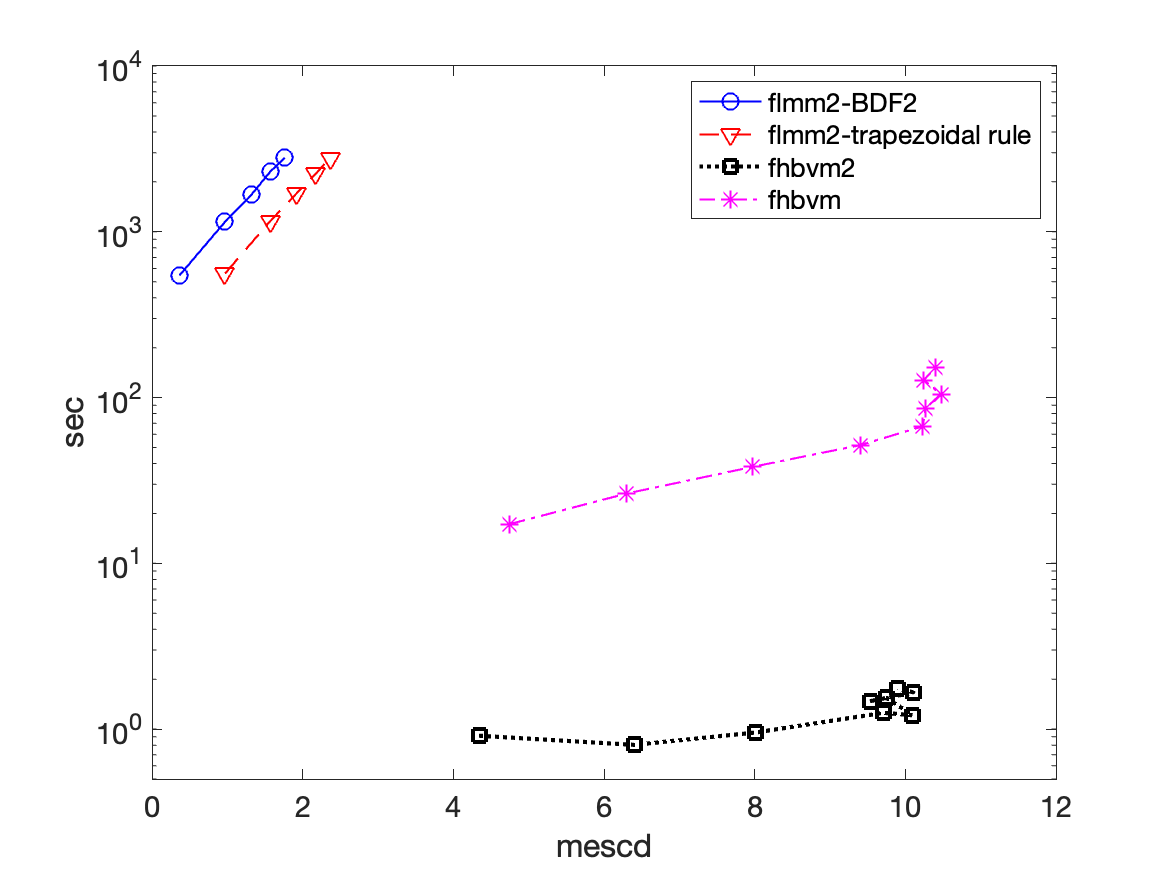}
\caption{Problem (\ref{ex3})-(\ref{ex3_Ay0}), WPD.}
\label{ex3wpd}
\end{figure}

Next, we consider the following stiff oscillatory problem:
\begin{equation}\label{ex3}
y^{(0.5)} = Ay,\qquad t\in[0,20],\qquad y(0) = y_0,
\end{equation}
with
\begin{equation}\label{ex3_Ay0}
A=\frac{1}{8}\pmatrix{rrrrr}
    41&   41&   -38&    40&    -2\\
   -79&    81&     2&     0&    -2\\
    20&  -60&   20&   -20&    -8\\
   -22&    58&   -24&    20&    -4\\
     1 &    1&    -2&    -4&    -2\\
         \endpmatrix, \qquad y_0= \pmatrix{c} 1\\ 2\\ 3\\ 4\\ 5\endpmatrix.
 \end{equation}    
The eigenvalues of $A$ are given by $\lambda_{1/2} =10\pm 10\ii$, $\lambda_{3/4}=\frac{1}2\pm\frac{1}2\ii$ and $\lambda_5=-1$: the first 4 eigenvalues are on the boundary of the stability region, $\partial\Lambda_{0.5}$, whereas $\lambda_5\in\Lambda_{0.5}$. The components of the solution, given by
$$y(t) =E_{0.5}(A\/t^{0.5})y_0, \qquad t\ge0,$$
are depicted in Figure~\ref{ex3fig},  with a fast oscillatory component superimposed on a slowly oscillating one. Moreover, it is clear that the vector field, alike the solution, is nonsmooth at $t=0$, so that a graded mesh is mandatory, in this case. We compare the codes with the following parameters:
\begin{itemize}
\item {\tt flmm2} (both when using the trapezoidal rule and the BDF2 method): 
$$h=\frac{20}{5\ell\cdot 10^5}, \qquad \ell=1,2,3,4,5;$$
\item {\tt fhbvm}: $M=50\ell$,\, $\ell=4,\dots,12$; 
\item {\tt fhbvm2}: $\nu=20$, $n=1$, $N=50\ell$,\, $\ell=4,\dots,12$. 
\end{itemize}
The obtained results are summarized in the WPD in Figure~\ref{ex3wpd}, which allows us to conclude that:
\begin{itemize}
\item the code {\tt flmm2} using the BDF2 method is the less effective one, able to reach less than 2 mescd in about  46 min;
\item the code {\tt flmm2} using the trapezoidal rule is able to reach about 2.4 mescd in about the same time;
\item the code {\tt fhbvm} is able to reach 10 mescd in about 67 sec (further reductions of the stepsize do not improve the accuracy, but only affect the execution time). Consequently, it is much more effective than the code {\tt flmm2};
\item the code {\tt fhbvm2} is able to obtain the same accuracy as that of the code {\tt fhbvm}, but in less than  2 sec (also now, further reductions of the stepsize do not improve the accuracy, but only affect the execution time).
\end{itemize}
From the above results, one concludes that the code {\tt fhbvm2} here introduced is the most effective. This is, indeed, expected because of the oscillating nature of the solution, for which the mixed-mesh strategy is very suited. 

\subsection{Problem~3}\label{p3}
\begin{figure}[ph]
\centering
\includegraphics[width=8cm]{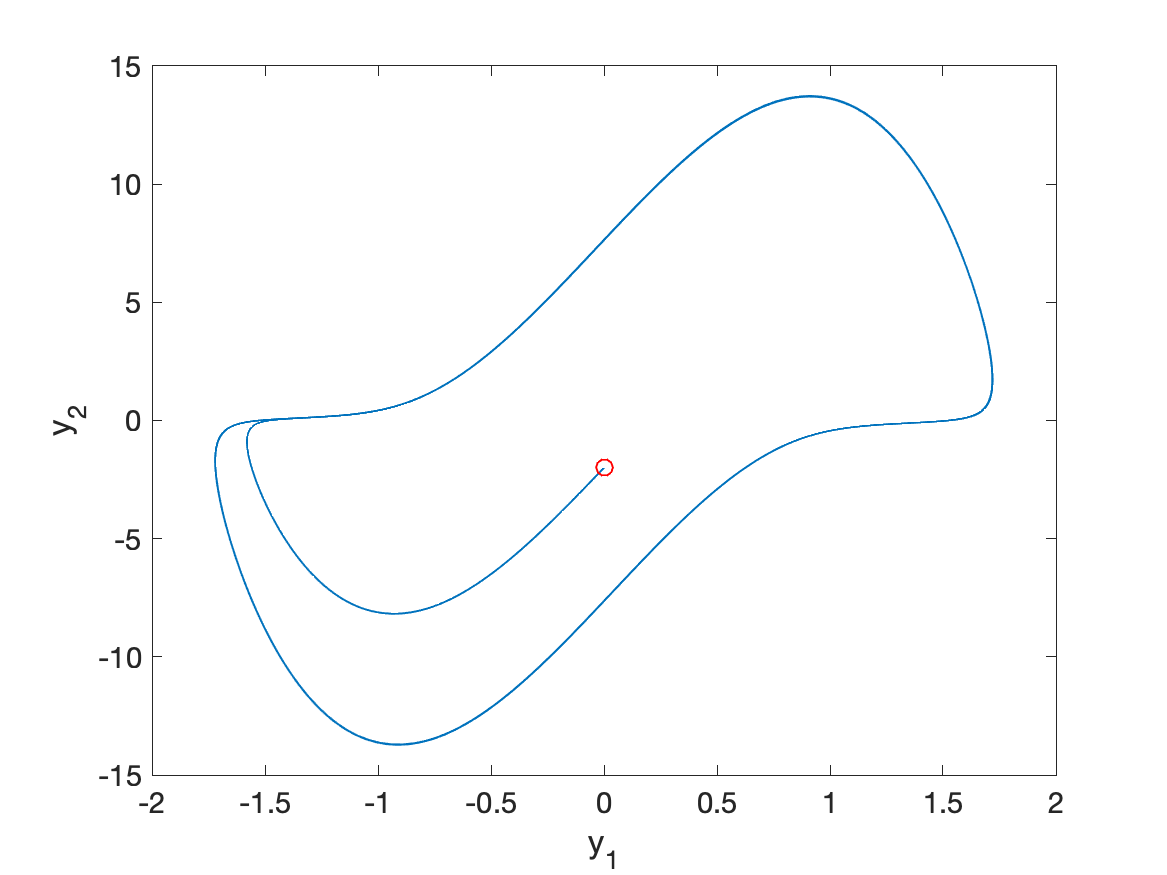}\includegraphics[width=8cm]{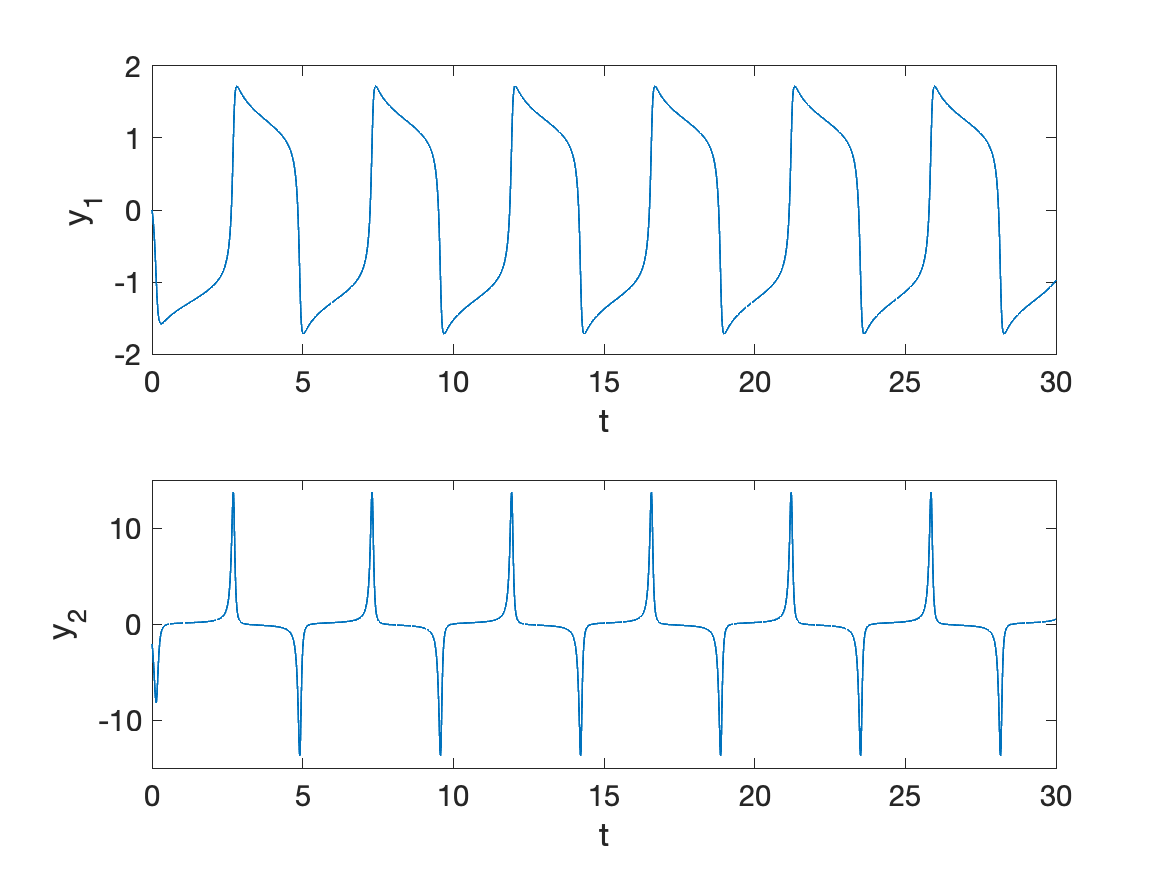}
\caption{Problem (\ref{ex2}),  phase portrait of the solution (left-plot) and solution versus time (right-plot). In the left-plot, the red circle is the starting point of the trajectory.}
\label{ex2fig}

\bigskip\bigskip
\includegraphics[width=10cm]{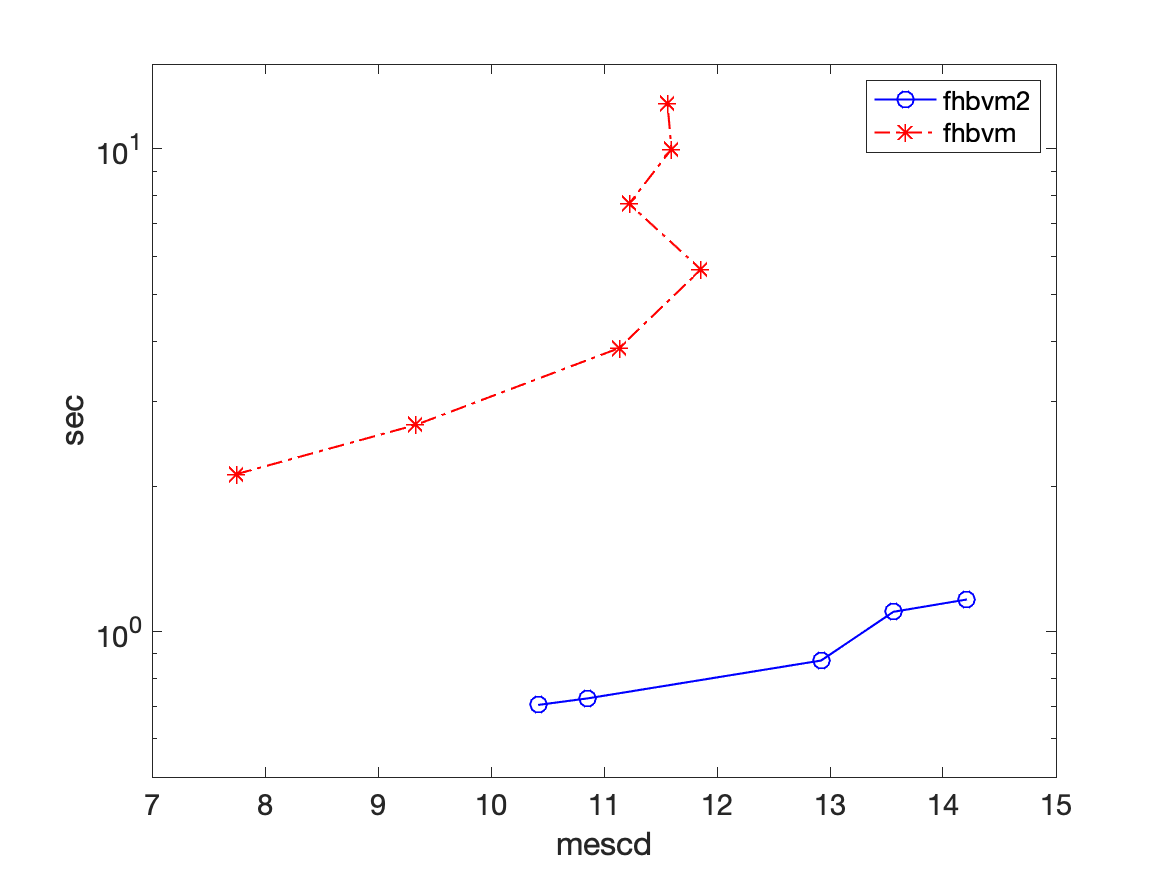}
\caption{Problem (\ref{ex2}), WPD.}
\label{ex2wpd}
\end{figure}

Next, we consider the following fractional version of the Van der Pol problem \cite{Petras2011}:
\begin{equation}\label{ex2}
y_1^{(0.9)} = y_2, \qquad y_2^{(0.9)} = -y_1-10y_2(y_1^2-1), \qquad t\in[0,30], \qquad y(0) = (0,\,-2)^\top.
\end{equation}
As in the ODE case, the solution approaches a limit cycle, as is shown in the left-plot of Figure~\ref{ex2fig}, so that the solution becomes eventually periodic (see also the right-plot in the same figure). For this problem, we do not consider comparisons with the code {\tt flmm2}, since it does not converge within a ``practical'' execution time. Instead, we compare the codes {\tt fhbvm} and {\tt fhbvm2}, by using the following parameters:\footnote{A reference solution at $t=30$ has been computed by using the code {\tt fhbvm2} with parameters $\nu=100$, $n=5$, $N=10^4$.} 
\begin{itemize}
\item {\tt fhbvm}: $M=50\ell$,\, $\ell=3,\dots,9$;
\item {\tt fhbvm2}:  $\nu=50$, $n=2$, $N=50\ell$,\, $\ell=4,\dots,8$.
\end{itemize}
The obtained results are depicted in the WPD of Figure~\ref{ex2wpd}, showing that the new code {\tt fhbvm2} can reach a higher accuracy, w.r.t. to the code {\tt fhbvm}, and with a smaller execution time. Also in this case, the mixed-mesh strategy implemented in the code {\tt fhbvm2} allows to efficiently cope with the periodic nature of the solution. 

\subsection{Problem~4}\label{p4}
\begin{figure}[ph]
\centering
\includegraphics[width=8cm]{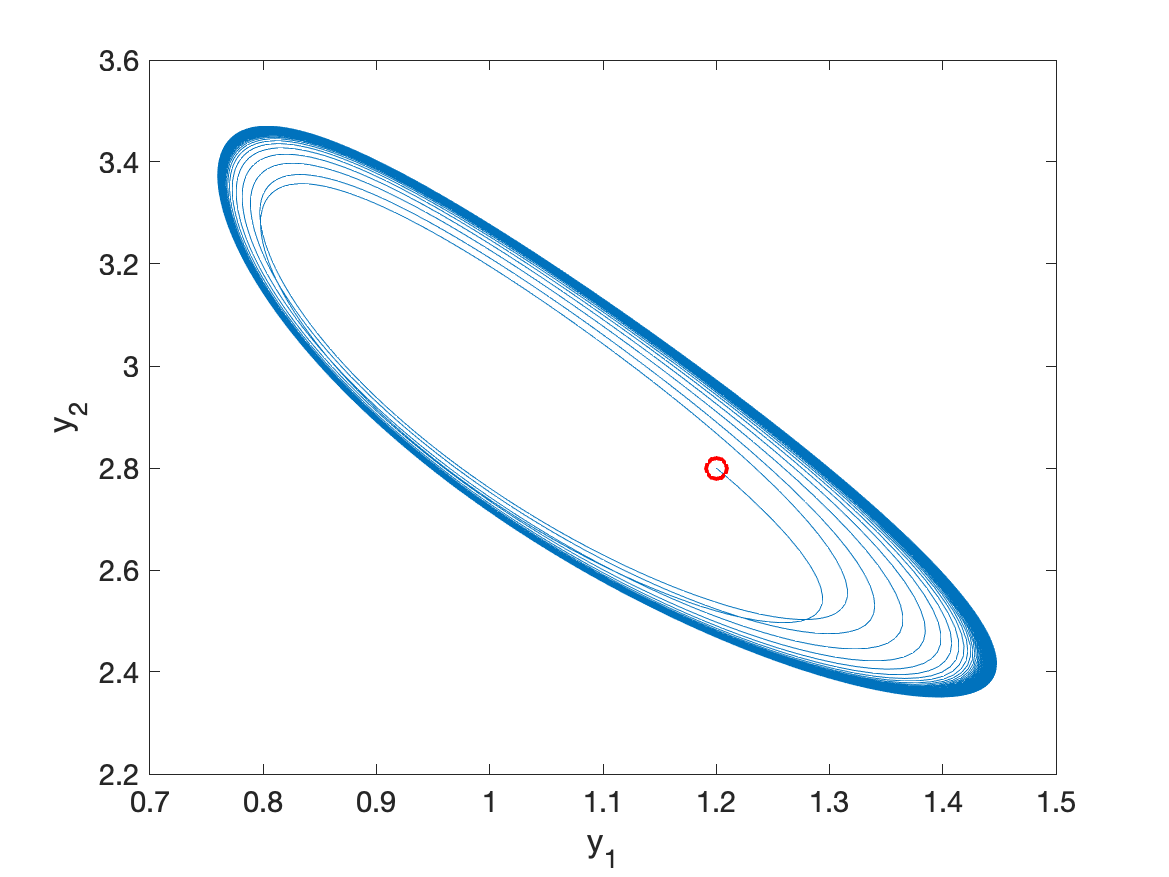}\includegraphics[width=8cm]{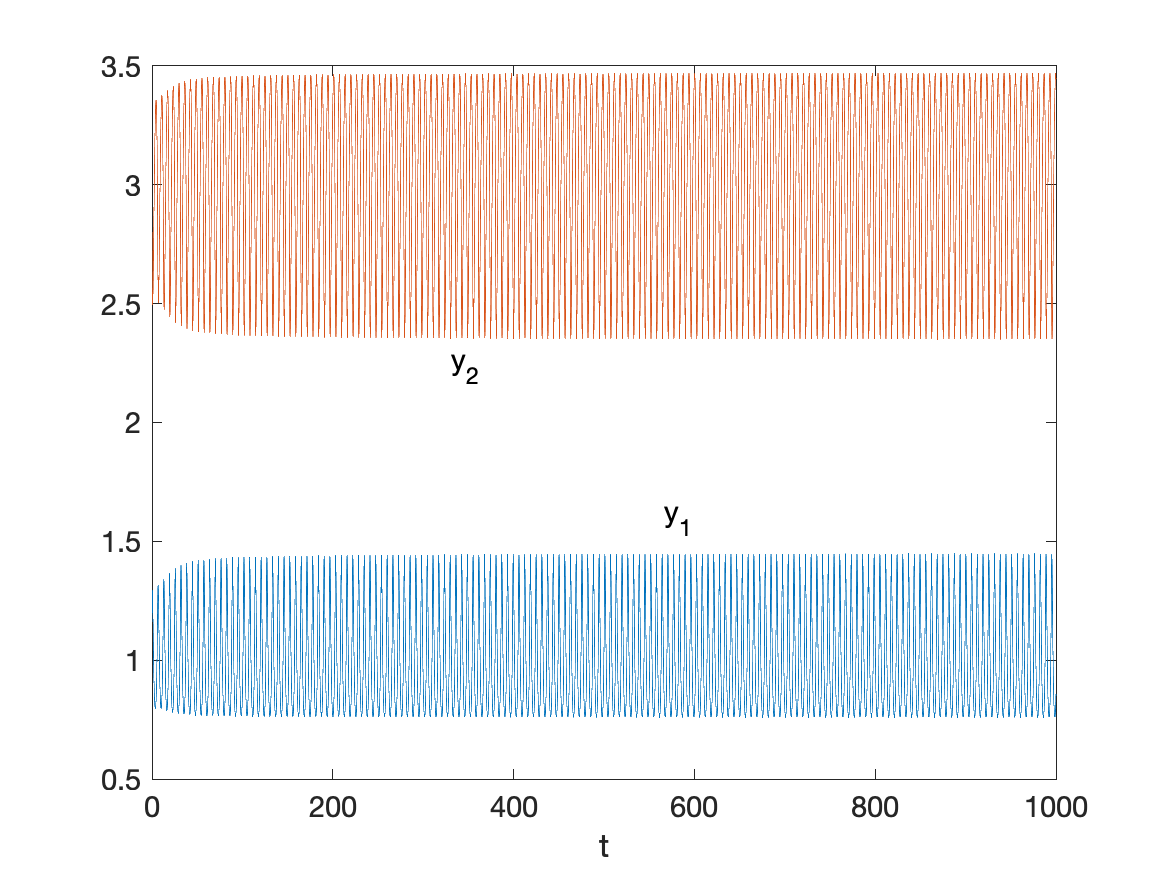}
\caption{Problem (\ref{ex4}),  phase portrait of the solution (left-plot) and solution versus time (right-plot). In the left-plot, the red circle is the starting point of the trajectory.}
\label{ex4fig}

\bigskip\bigskip
\includegraphics[width=10cm]{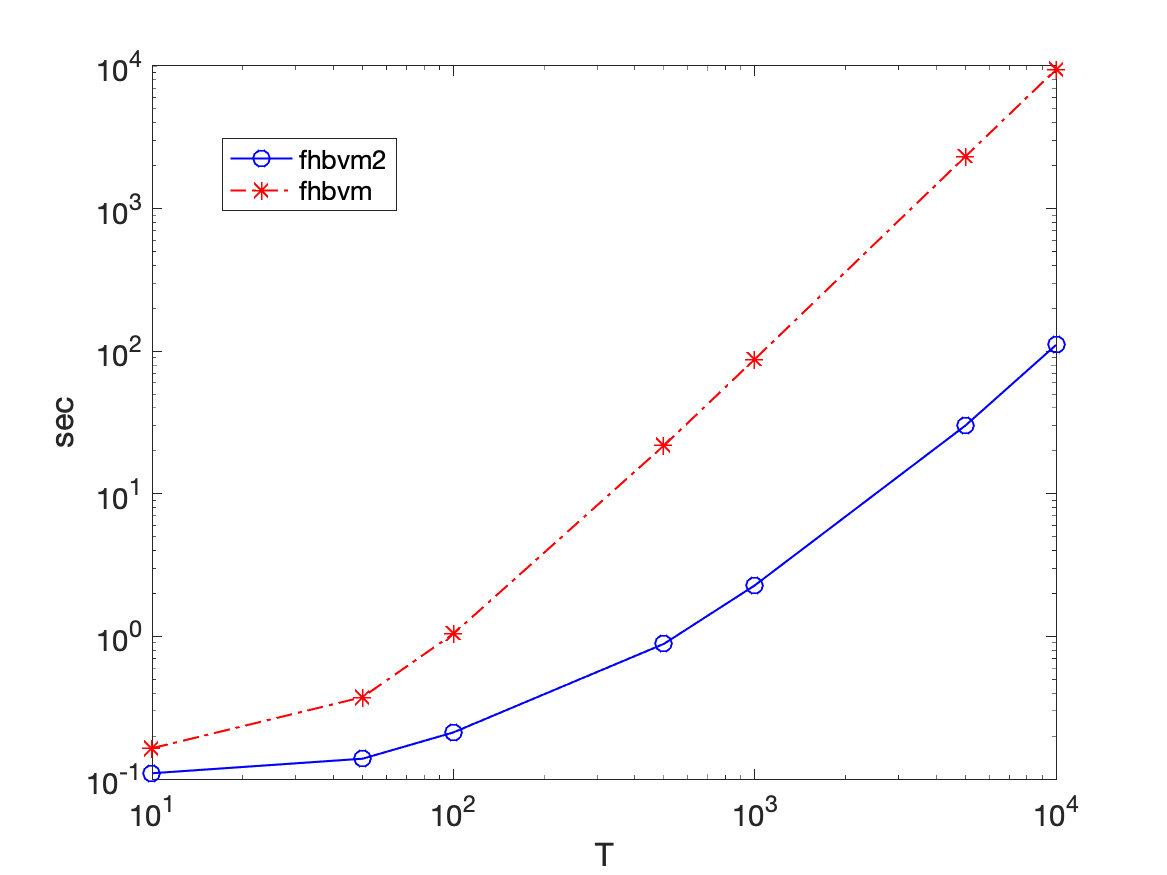}
\caption{Problem (\ref{ex4}),  execution time versus $T$.}
\label{ex4fig1}
\end{figure}

At last, we consider the following fractional Brusselator problem \cite{BGI2024}:
\begin{equation}\label{ex4}
y_1^{(0.7)} =   1-4y_1+y_1^2 y_2, \qquad y_2^{(0.7)} = 3y_1-y_1^2y_2,\qquad t\in[0,T],\qquad y(0)=(1.2,~2.8)^\top,
\end{equation}
for which we select increasing widths $T$ of the time interval ($T$ is chosen as an integer value). For this problem, the solution approaches a limit cycle, so that it is of periodic type, as is shown in the two plots in Figure~\ref{ex4fig} for $T=1000$.

Our aim is that of comparing the code {\tt fhbvm2}, used with parameters $\nu=20, n=1, N=T$ (so that the timestep in the uniform mesh is $h=1$) with the code {\tt fhbvm}, used with parameter $M=T$ (so that the final timestep used is approximately equal to 1), for increasing values of $T$.\footnote{Also, for this problem, we do not report comparisons with the code {\tt flmm2}, since to obtain comparable accuracies w.r.t. the other codes, too high execution times would be needed.}
 
As expected, the two codes compute, up to a level compatible with round-off errors (which, however, increase with $T$), the same approximation to $y(T)$, but with different execution times, which we want to compare: they are plot in Figure~\ref{ex4fig1}, showing that {\tt fhbvm2} outperforms {\tt fhbvm}.  Moreover, in Table~\ref{ex4tab} we list the number of mesh-points, needed by the two codes, w.r.t. $T$: one verifies that for larger values of $T$, the number of mesh-points used by {\tt fhbvm} is roughly 10 times larger than those used by the new code {\tt fhbvm2}. As is expected, also in this case, the usefulness of the mixed-mesh strategy, implemented in the code {\tt fhbvm2}, turns out to be very effective, to cope with the periodic nature of the solution.

\begin{table}[t]
\centering
\caption{Problem (\ref{ex4}), number of mesh-points required by the codes for increasing values of $T$.}
\label{ex4tab}\smallskip
\begin{tabular}{|r|r|r|}
\hline
$T$ & {\tt fhbvm} & {\tt fhbvm2}\\
\hline
$10$ & 95& 30\\
$50$ & 483 & 70\\
$100$ & 968 & 120 \\
$500$ & 4850 &520 \\
$1000$ & 9702& 1020\\
$5000$ & 48521 & 5020\\
$10000$ &97044 & 10020\\
\hline
\end{tabular}
\end{table}

\section{Conclusions}\label{fine}
In this paper we have provided a linear stability analysis of FHBVM methods  solving Caputo FDE-IVPs, thus confirming the effectiveness of the higher-order methods in approximating the one parameter Mittag-Leffler function. Moreover, we have given the implementation details of the methods when using a mixed-mesh strategy, consisting in coupling an initial graded mesh with a subsequent uniform one. This mixed-mesh implementation, in turn, is aimed at effectively coping with problems, having a nonsmooth vector field at the origin and with solutions of oscillatory type, over wide time-intervals. Numerical tests show the usefulness of this strategy, implemented in the Matlab\cpr code {\tt fhbvm2}, especially when solving FDE problems with periodic solutions. In this respect, the new code {\tt fhbvm2} represents a noticeable improvement over the previous one, {\tt fhbvm}: the improvement is even more impressive, if compared with other existing available codes.

\paragraph*{Data availability.} The code {\tt fhbvm2} is available at the URL \cite{fhbvm}.
 
\paragraph*{\bf Acknowledgements.} The first three authors are members of the ``Gruppo Nazionale per il Calcolo Scientifico-Istituto Nazionale di Alta Matematica (GNCS-INdAM)''. The last author is supported by the project n.ro PUTJD1275 of the Estonian Research Council.

\paragraph*{Declarations.} The authors declare no conflict of interests.

\end{document}